\documentclass[12pt]{amsart}
\usepackage{amssymb}

\theoremstyle{plain}
\newtheorem{theorem}{Theorem}[section]
\newtheorem{corollary}[theorem]{Corollary}
\newtheorem{lemma}[theorem]{Lemma}
\newtheorem{proposition}[theorem]{Proposition}
\newtheorem*{question}{Question}

\theoremstyle{definition}
\newtheorem{definition}[theorem]{Definition}
\newtheorem{remark}[theorem]{Remark}
\newtheorem{numbered}[theorem]{}
\newtheorem{example}[theorem]{Example}

\font\bi=cmbxti12 % /usr/local/teTeX/texmf/fonts/source/public/cm

\def\mid{\::\:} % semicolon as a separator
\def\tnun{{\textstyle{\frac{T^n}{\nu^n}}}}  % a typical fraction
\def\su{{\scriptscriptstyle U}} % often used superscript

\newcommand{\marg}[1]{}
\newcommand{\lab}[1]{\label{#1}\marg{#1}}
\newcommand{\abs}[1]{\lvert#1\rvert}
\newcommand{\norm}[1]{\lVert#1\rVert}
\newcommand{\bigabs}[1]{\bigl\lvert#1\bigr\rvert}

\newcommand{\Bignorm}[1]{\Bigl\lVert#1\Bigr\rVert}
\renewcommand{\le}{\leqslant}
\renewcommand{\ge}{\geqslant}
\newcommand{\msn}[2]{{\mathfrak m}_{#1#2}} % mixed seminorm
\newcommand{\term}[1]{{\bi #1}}
\newcommand{\ssection}[1]{\medskip\leftline{\bf #1}
                          \nopagebreak\smallskip\nopagebreak}
\DeclareMathOperator{\Null}{Null}

\title[Spectral radii]
      {Spectral radii of bounded operators\\
       on topological vector spaces}
\author{Vladimir G. Troitsky}
\address{Department of Mathematics, C1200
         The University of Texas
         Austin, TX 78712. USA}
\email{vladimir@math.utexas.edu}
\date{\today}

\begin{document}

\baselineskip 16pt

\begin{abstract}
  In this paper we develop a version of spectral theory 
  for bounded linear operators on topological vector spaces. We show
  that the Gelfand formula for spectral radius and Neumann series can
  still be naturally interpreted for operators on topological vector
  spaces. Of course, the resulting theory has many similarities to the 
  conventional spectral theory of bounded operators on Banach spaces,
  though there are several important differences. The main difference
  is that an operator on a topological vector space has several
  spectra and several spectral radii, which fit a well-organized
  pattern.
\end{abstract}

\maketitle
%\tableofcontents
\setcounter{section}{-1}

\section{Introduction}

The spectral radius of a bounded linear operator $T$ on a Banach space
is defined by the Gelfand formula $r(T)=\lim_n\sqrt[n]{\norm{T^n}}$.
It is well known that $r(T)$ equals the actual radius of the spectrum
$\bigabs{\sigma(T)}=\sup\{\abs{\lambda}\mid\lambda\in\sigma(T)\}$.
Further, it is known that the resolvent $R_\lambda=(\lambda I-T)^{-1}$
is given by the Neumann series $\sum_{i=0}^\infty
\frac{T^i}{\lambda^{i+1}}$ whenever $|\lambda|>r(T)$. It is natural to
ask if similar results are valid in a more general setting, e.g., for
a bounded linear operator on an arbitrary topological vector space.
The author arrived to these questions when generalizing some results
on Invariant Subspace Problem from Banach lattices to ordered
topological vector spaces.  One major difficulty is that it is not
clear which class of operators should be considered, because there are
several non-equivalent ways of defining bounded operators on
topological vector spaces. Another major difficulty is the lack of a
readily available developed spectral theory. The spectral theory of
operators on Banach spaces has been thoroughly studied for a long
time, and is extensively used.  Unfortunately, little has been known
about spectral theory of bounded operators on general topological
vector spaces, and many techniques used in Banach spaces cannot be
applied for operators on topological vector spaces.  In particular,
the spectrum, the spectral radius, and the Neumann series are the
tools which are widely used in the study of the Invariant Subspace
Problem in Banach spaces, but which have not been sufficiently studied
for general topological vector spaces. To overcome this obstacle we
have developed a version of the spectral theory of bounded operators
on general topological vector spaces and on locally convex spaces.
Some results in this direction have also been obtained by B.~Gramsch
\cite{Gramsch:66}, and by F.~Garibay and R.~Vera
\cite{Garibay:97,Garibay:98,Vera:97}.

In particular, we consider the following classification of bounded
operators on a topological vector space. We call a linear operator~$T$
\begin{itemize}
  \item[-] \term{nb-bounded} if $T$ maps some neighborhood of zero into
    a bounded set,
  \item[-] \term{nn-bounded} if there
    is a base of neighborhoods of zero such that $T$ maps
    every neighborhood in this base into a multiple of itself, and
  \item[-] \term{bb-bounded} if $T$ maps bounded sets into bounded sets.  
\end{itemize}
The classes of all linear operators, of all bb-bounded operators, of
all continuous operators, of all nn-bounded operators, and of all
nb-bounded operators form nested algebras. The spectrum of an operator
$T$ in each of these algebras is defined as usual, i.e., the set of
$\lambda$'s for which $\lambda I-T$ is not invertible in this algebra.
We show that the well known Gelfand formula for the spectral radius of
an operator on a Banach space
$r(T)=\lim\limits_{n\to\infty}\sqrt[n]{\norm{T^n}}$ can be generalized
to each of the five classes of operators on topological vector spaces,
and then we use this formula to define the spectral radius of an
operator in each of the classes. Then in Section~\ref{s:s-sr} we show
that if $T$ is a continuous operator on a sequentially complete
locally convex space and $\abs{\lambda}$ is greater than the spectral
radius of $T$ in any of the five classes, then the Neumann series
$\sum_{n=0}^{\infty}\frac{T^n}{\lambda^{n+1}}$ converges in the
topology of the class, and $\lambda$ does not belong to the
corresponding spectrum of $T$, i.e., the spectral radius is greater
than or equal to the geometrical radius of the spectrum. In
Sections~\ref{s:nb-bdd} and~\ref{s:compact} we show that the radii are
equal for nb-bounded and compact operators.
  
This paper is based on a part of the author's Ph.D.
thesis~\cite{Troitsky:PhD}. I would like to thank Yuri Abramovich for
his encouragement and for numerous suggestions and improvements.
Robert Kaufman read parts of the draft and asked interesting
questions, which resulted in new developments. Thanks are also due to
Heinrich Lotz, Michael Neumann, and Joseph Rosenblatt for their support
and interest in my work.

\section{Preliminaries and notation}  \lab{s:notation}

The symbols $X$ and $Y$ always denote topological vector spaces.  A
\term{neighborhood} of a point $x\in X$ is any subset of $X$
containing an open set which contains $x$. Neighborhoods of zero will
often be referred to as \term{zero neighborhoods}. Every zero
neighborhood $V$ is \term{absorbing}, i.e., $\bigcup_{n=1}^\infty
nV=X$.  In every topological vector space (over $\mathbb R$ or
$\mathbb C$) there exists a base ${\mathcal N}_0$\lab{a:N0} of zero
neighborhoods with the following properties:
  \begin{enumerate}
    \item Every $V\in{\mathcal N}_0$ is \term{balanced}, i.e., $\lambda
      V\subseteq V$ whenever $|\lambda|\le 1$; 
    \item For every $V_1,V_2\in{\mathcal N}_0$ there exists
      $V\in{\mathcal N}_0$ such that $V\subseteq V_1\cap V_2$;
    \item For every $V\in{\mathcal N}_0$ there exists $U\in{\mathcal N}_0$
      such that $U+U\subseteq V$;
    \item For every $V\in{\mathcal N}_0$ and every scalar $\lambda$ the
      set $\lambda V$ is in ${\mathcal N}_0$.
  \end{enumerate}
Whenever we mention a base zero neighborhood, we
assume that the base satisfies these properties. 

A topological vector space is called \term{normed} if the topology is
given by a norm. In this case the collection of all balls centered at
zero is a base of zero neighborhoods. A complete
normed space is referred to as a \term{Banach space}.
See~\cite{Dunford:58} for a detailed study of normed and Banach
spaces.

A subset $A$ of a topological vector space is called \term{bounded} if
it is absorbed by every zero neighborhood, i.e., for every zero
neighborhood $V$ one can find $\alpha>0$ such that $A\subseteq\alpha
V$. A set $A$ in a topological vector space is said to be
\term{pseudo-convex} or \term{semi-convex} if $A+A\subseteq\alpha A$
for some number $\alpha$ (for convex sets $\alpha=2$). If $U$ is a
zero neighborhood, $(x_\gamma)$ is a net in $X$, and $x\in X$, we
write $x_\gamma\xrightarrow{\su}x$\lab{a:U-conv} if for every $\varepsilon>0$
one can find an index $\gamma_0$ such that $x_\gamma-x\in\varepsilon U$
whenever $\gamma\ge\gamma_0$. It is easy to see that when $U$ is
pseudo-convex, this convergence determines a topology on $X$, and the
set of all scalar multiples of $U$ forms a base of the topology.
We denote $X$ equipped with this topology by $(X,U)$. 
Clearly, $(X,U)$ is Hausdorff if and
only if $\bigcap_{n=1}^\infty\frac{1}{n}U=\{0\}$.

A topological vector space is said to be \term{locally bounded} if
there exists a bounded zero neighborhood. Notice that if $U$ is a
bounded zero neighborhood then it is pseudo-convex. Conversely, if $U$
is a pseudo-convex zero neighborhood, then $(X,U)$ is locally bounded.
Recall that a \term{quasinorm} is a real-valued function on a vector
space which satisfies all the axioms of norm except the triangle
inequality, which is substituted by $\norm{x+y}\le
k\,\bigl(\norm{x}+\norm{y}\bigr)$ for some fixed positive constant
$k$. It is known (see, e.g.,~\cite{Kothe:60}) that a
topological vector space is quasinormable if and only if it is locally
bounded and Hausdorff. A complete quasinormed space is called
\term{quasi-Banach}.

If the topology of a topological vector space $X$ is given by a
\term{seminorm} p, we say that $X=(X,p)$ is a \term{seminormed}
space. Clearly, in this case $X=(X,U)$ where the convex set $U$ is the
unit ball of $p$ and, conversely, $p$ is the Minkowski functional of
$U$. A Hausdorff topological vector space is called \term{locally
convex} if there is a base of convex zero neighborhoods or,
equivalently, if the topology is generated by a family of seminorms
(the Minkowski functionals of the convex zero neighborhoods).  When
dealing with locally convex spaces we will always assume that the base
zero neighborhoods are convex. Similarly, a Hausdorff topological
vector space is said to be \term{locally pseudo-convex} if it has a
base of pseudo-convex zero neighborhoods.  A complete metrizable
topological vector space is usually referred to as a \term{Fr\'echet
space}.

Further details on topological vector spaces can be found in~
\cite{Dunford:58,Kothe:60,Robertson:64,Edwards:65,Schaefer:71,Kelley:76}.
For details on locally bounded
and quasinormed topological vector spaces we refer the reader
to~\cite{Kothe:60,Kalton:84,Rolewicz:85}.

By an \term{operator} we always mean a \term{linear operator} between
vector spaces. We will usually use the symbols $S$ and $T$ to denote
operators. Recall that an operator $T$ between normed spaces is said
to be \term{bounded} if its operator norm defined by
$\norm{T}=\sup\{\norm{Tx}\mid \norm{x}\le 1\}$ is finite. It is well
known that an operator between normed spaces is bounded if and only if
it is continuous. An operator between two vector spaces is said to be
of \term{finite rank} if the range of $T$ is finite dimensional.

If $\mathcal A$ is a unital algebra and $a\in\mathcal A$, then the
\term{resolvent set}\lab{a:rho} of $a$ is the set $\rho(a)$ of all
$\lambda\in\mathbb C$ such that $e-\lambda a$ is invertible in
$\mathcal A$. The resolvent set of an element $a$ in a non-unital
algebra $\mathcal A$ is defined as the set of all $\lambda\in\mathbb
C$ for which $e-\lambda a$ is invertible in the unitalization
${\mathcal A}_\times$ of $\mathcal A$.  The
\term{spectrum}\lab{a:spectrum} of an element of an algebra is defined
via $\sigma(a)={\mathbb C}\setminus\rho(a)$. It is well-known that
whenever $\mathcal A$ is a unital Banach algebra then $\sigma(a)$ is
compact and nonempty for every $a\in\mathcal A$. In this case the
\term{spectral radius}\lab{a:sr} $r(a)$ is defined via \term{Gelfand
  formula}: $r(a)=\lim\limits_{n\to\infty}\sqrt[n]{\norm{a^n}}$. It
is well known that $r(a)=\bigabs{\sigma(a)}$, where
$\bigabs{\sigma(a)}$ is the \term{geometrical radius}\lab{a:grsp} of
$\sigma(a)$, i.e.,
$\bigabs{\sigma(a)}=\sup\{\abs{\lambda}\mid\lambda\in\sigma(a)\}$. An
element $a\in A$ is said to be \term{quasinilpotent} if
$\sigma(a)=\{0\}$.

If $T$ is a bounded operator on a Banach space $X$ then we will
consider the spectrum $\sigma(T)$ and the resolvent set $\rho(T)$ in
the sense of the Banach algebra of bounded operators on $X$. If
$\lambda\in\rho(T)$ then the inverse $(I-\lambda T)^{-1}$ is called
the \term{resolvent operator}\lab{a:res} and is denoted by
$R(T;\lambda)$ or just $R_\lambda$. It is well known that if
$\lambda\in\mathbb C$ satisfies $\abs{\lambda}>r(T)$ then the
\term{Neumann series} $\sum_{i=0}^\infty \frac{T^i}{\lambda^{i+1}}$
converges to $R_\lambda$ in operator norm.  We say that $T$ is
\term{locally quasinilpotent} at $x\in X$ if
$\lim\limits_{n\to\infty}\sqrt[n]{\norm{T^nx}}=0$.

\section{Bounded operators} \lab{s:bounded-ops}

There are various definitions for a bounded linear operator between
two topological vector spaces. To avoid confusion, we will, of course,
give different names to different types of boundedness.

\begin{definition}
  Let $X$ and $Y$ be topological vector spaces. An operator $T\colon
  X\to Y$ is said to be
  \begin{enumerate}
    \item \term{bb-bounded} if it maps every
    bounded set into a bounded set; \lab{i:def-bb}
   \item  \term{nb-bounded} if it maps some neighborhood into
      a bounded set; \lab{i:def-nb}
  \end{enumerate}
  Further, if $X=Y$ we will say that $T\colon X\to X$ is
  \term{nn-bounded} if there exists a base ${\mathcal N}_0$ of zero
  neighborhoods such that for every $U\in{\mathcal N}_0$
  there is a positive scalar $\alpha$ such that
  $T(U)\subseteq \alpha U$.
 \end{definition}

\begin{remark}
  \cite{Edwards:65} and~\cite{Kelley:76} present (\ref{i:def-bb}) as the
  definition of a bounded operator on a topological vector space,
  while~\cite{Robertson:64} and~\cite{Schaefer:71} use~(\ref{i:def-nb}) for the
  same purpose. As we will see, these definitions are far from being
  equivalent.
\end{remark}

\begin{proposition} \lab{p:hierarchy}
  Let $X$ and $Y$ be topological vector spaces. For an operator
  $T\colon X\to Y$ consider the following statements:
  \begin{enumerate}
    \item $T$ is bb-bounded; \lab{pi:hierarchy:bb}
    \item $T$ is continuous; \lab{pi:hierarchy:cont}
    \item $T$ is nn-bounded; \lab{pi:hierarchy:nn}
    \item $T$ is nb-bounded. \lab{pi:hierarchy:nb}
  \end{enumerate}
  Then
  {\rm(\ref{pi:hierarchy:nb}) $\Rightarrow$
      (\ref{pi:hierarchy:cont}) $\Rightarrow$
      (\ref{pi:hierarchy:bb})}. 
  Furthermore, if $X=Y$ then
  {\rm(\ref{pi:hierarchy:nb}) $\Rightarrow$
      (\ref{pi:hierarchy:nn}) $\Rightarrow$
      (\ref{pi:hierarchy:cont}) $\Rightarrow$
      (\ref{pi:hierarchy:bb})}.
\end{proposition}

\begin{proof} 
%  If $T$ is nb-bounded, then $T(U)$ is bounded for some zero neighborhood
%  $U$ in $X$. Therefore for every zero neighborhood $V$ of $Y$ we have
%  $T(U)\subseteq\alpha V$ for some $\alpha>0$. Then $T(\alpha^{-1}U)
%  \subseteq V$, so that $T$ is continuous.
%
%  Assume now that $T$ is continuous, and let $A$ be a bounded subset
%  of $X$. If $V$ be a neighborhood of zero in
%  $Y$, then there exists a neighborhood of zero $W$ in $X$ such that
%  $T(W)\subseteq V$, but since $A\subseteq\alpha W$ for some positive
%  $\alpha$ then $T(A)\subseteq \alpha T(W)\subseteq\alpha V$, so that
%  $T$ is bb-bounded.
  The implications 
  {\rm(\ref{pi:hierarchy:nb}) $\Rightarrow$
      (\ref{pi:hierarchy:cont}) $\Rightarrow$
      (\ref{pi:hierarchy:bb})} 
  are trivial. To show 
  {\rm(\ref{pi:hierarchy:nb}) $\Rightarrow$
      (\ref{pi:hierarchy:nn}) $\Rightarrow$
      (\ref{pi:hierarchy:cont})}
  assume that $X=Y$ and fix a base ${\mathcal N}_0$ of zero
  neighborhoods. If $T$ is nb-bounded then $T(U)$ is bounded for some
  $U\in{\mathcal N}_0$. Note that $\widetilde{\mathcal N}_0=\{\,\lambda U\cap
  V\mid \lambda>0,~V\in{\mathcal N}_0\,\}$ is another base of zero
  neighborhoods. For each $W=\lambda U\cap V$ in 
  $\widetilde{\mathcal N}_0$ we have $T(W)\subseteq \lambda T(U)$.
  But $T(U)$ is
  bounded and so $T(W)\subseteq \lambda T(U)\subseteq \lambda\alpha
  W$ for some positive $\alpha$, i.e., $T$ is nn-bounded.%
  \footnote{Note that if the topology is locally convex, then we can
  assume that $U$ is convex and ${\mathcal N}_0$ consists of convex
  neighborhoods. In this case $\widetilde{\mathcal N}_0$ also consists of
  convex neighborhoods.}
%  $\widetilde{V}\in\widetilde{\mathcal N_0}$, then
%  $\widetilde{V}=nU\cap V$ for some
%  $n\in\mathbb N$ and $V\in\mathcal N_0$, so that there exists a
%  neighborhood $W$ in $\mathcal N_0$ such that $W\subseteq nU\cap
%  V=\widetilde{V}$. Since $T(U)$ is bounded we have
%  \begin{displaymath}
%    T(\widetilde{V})\subseteq T(nU\cap V)\subseteq nT(U)\subseteq n\alpha
%    W\subseteq n\alpha\widetilde{V}
%  \end{displaymath} 
%  for some positive $\alpha$, so that $T$ is nn-bounded.
 
  Finally, if $T$ is nn-bounded, then there is a base
  ${\mathcal N}_0$ such that for every zero neighborhood
  $U\in{\mathcal N}_0$ there is a positive scalar $\alpha$
  such that $T(U)\subseteq \alpha U$.  Let $V$ be an arbitrary zero
  neighborhood. Then there exists $U\in{\mathcal N}_0$ such
  that $U\subseteq V$, so that $T(U)\subseteq\alpha U\subseteq\alpha
  V$ for some $\alpha>0$. Taking $W=\frac{1}{\alpha}U$ we get
  $T(W)\subseteq V$, hence $T$ is continuous.
\end{proof}

\begin{numbered} \lab{n:lb-equal}
  It can be easily verified that if $T$ is an operator on a locally
  bounded space then all the statements in Lemma~\ref{p:hierarchy} are
  equivalent. In general, however, these notions are not equivalent.
  Obviously, the identity operator $I$ is always nn-bounded, continuous,
  and bb-bounded, but $I$ is
  nb-bounded if and only if the space is locally bounded.  Every
  bb-bounded operator between two locally convex spaces is continuous
  if and only if the domain space is bornological.  (Recall that a
  locally convex space is \term{bornological} if every balanced
  convex set absorbing every bounded set is a \marg{Every metrizable
    locally convex space is bornological.} zero neighborhood, for
  details see~\cite{Schaefer:71,Robertson:64}.)
\end{numbered}

\begin{example}{\it A continuous but not nn-bounded operator.}
  Let $T$ be the left shift on the space of all real sequences
  ${\mathbb R}^{\mathbb N}$ with the topology of coordinate-wise
  convergence, i.e.,
  $T\colon(x_0,x_1,x_2,\dots)\mapsto(x_1,x_2,x_3,\dots)$. Clearly $T$
  is continuous. We will show that $T$ is not nn-bounded. Assume that
  for every zero neighborhood $U$ in some base ${\mathcal N}_0$ there
  is a positive scalar $\alpha$ such that $T(U)\subseteq \alpha U$.
  Since the set $\bigl\{x=(x_k)\mid\abs{x_0}<1\bigr\}$ is a zero
  neighborhood, there must be a base neighborhood $U\in{\mathcal N}_0$
  such that $U\subseteq\{x\mid\abs{x_0}<1\}$. Since $T(U)\subseteq
  \alpha U$ for some positive $\alpha$ then $T^n(U)\subseteq \alpha^n
  U$, so that if $x=(x_k)\in U$ then $T^nx\in\nolinebreak\alpha^nU$, so that
  $\abs{x_n}=\bigabs{(T^nx)_0}<\alpha^n$. Hence
  $U\subseteq\bigl\{x\mid \abs{x_n}<\alpha^n\mbox{ for each }
  n>0\bigr\}$. But this set is bounded, while the space
  is not locally bounded, a contradiction.
\end{example}

%  In order to construct an example of an nn-bounded operator which is
%  not nb-bounded, consider the space $C(0,1)$ of continuous functions
%  on the unit interval with pointwise convergence topology. It is a
%  locally convex topology with $0$-neighborhoods of the form $\{\,f\in
%  C(0,1)\mid |f(x_i)|<\varepsilon; x_1,\dots,x_n\in[0,1]\,\}$.
%  Consider the operator on $C(0,1)$ that maps $f\in C(0,1)$ into the
%  constant function $f(0)$. This operator is nb-bounded but not
%  nn-bounded.

\begin{numbered} \lab{n:algebra}
  {\bf Algebraic properties of bounded operators}. The sum of two
  bb-bounded operators is bb-bounded because the sum of two bounded
  sets in a topological vector space is bounded. Clearly the
  product of two bb-bounded operators is bb-bounded. It is well known
  that sums and products of continuous operators are continuous.
  Obviously, the product of two nn-bounded operators is nn-bounded,
  and it can be easily verified that the sum of two nn-bounded
  operators on a locally convex (or locally pseudo-convex) space is
  again nn-bounded. It is not difficult to see that the sum of two
  nb-bounded operator is nb-bounded. Indeed, suppose that $T_1$ and
  $T_2$ are two nb-bounded operators, then the sets $T_1(U_1)$ and
  $T_2(U_2)$ are bounded for some base zero neighborhoods $U_1$ and
  $U_2$. There exists another base zero neighborhood $U\subseteq
  U_1\cap U_2$, then the sets $T_1(U)$ and $T_2(U)$ are bounded, so
  that $(T_1+T_2)(U)\subseteq T_1(U)+T_2(U)$ is bounded. Finally, it
  is not difficult to see that the product of two nb-bounded operators
  is again nb-bounded. In fact, it follows immediately from
  Proposition~\ref{p:hierarchy} and the following simple observation:
  if we multiply an nb-bounded operator by a bb-bounded operator on
  the left or by an nn-bounded operator on the right, the product is
  nb-bounded.
  
  Thus, the class of all bb-bounded operators, the class of all
  continuous operators, and the class of all nb-bounded operators are
  subalgebras of the algebra of all linear operators. The class of
  nn-bounded operators is an algebra provided that the space is
  locally (pseudo-)convex.
\end{numbered}

\ssection{Boundedness in terms of convergence}
  Suppose $T\colon X\to Y$ is an operator between two topological
  vector spaces. It is well known that $T$ is continuous if and only
  if it maps convergent nets to convergent nets.

  Notice that a subset of a topological vector space is unbounded if and
  only if it contains an unbounded sequence. Therefore, an operator is
  bb-bounded if and only if it maps bounded sequences (nets) to
  bounded sequences (respectively nets).
  
%  It is not difficult to see that $T$ is bb-bounded if and only if it
%  maps bounded nets to bounded nets. Indeed, if $T$ is bb-bounded, and
%  $(x_\alpha)$ is a bounded net, then $\{x_\alpha\}$ is bounded as a
%  set, so that the set $\{Tx_\alpha\}$ is bounded hence the net
%  $(Tx_\alpha)$ is bounded.  Reversely, suppose $T$ maps bounded nets
%  to bounded nets and assume \marg{Exercise for the reader?} that
%  $T(A)$ is unbounded for some bounded set $A$. Then there is a zero
%  neighborhood $V$ in $Y$ such that $V$ does not absorb $T(A)$, that
%  is for every $n\ge 1$ there exists $y_n\in T(A)$ so that $y_n\notin
%  nV$. Suppose $y_n=Tx_n$ for some $x_n\in A$ for every $n$. Then $T$
%  maps a bounded sequence $(x_n)$ to an unbounded sequence $(y_n)$, a
%  contradiction.
  
  It is easy to see that $T$ is nn-bounded if and only if $T$ maps
  $U$-bounded ($U$-convergent to zero) sequences to $U$-bounded (respectively
  $U$-convergent to zero) sequences for every base zero neighborhood $U$ in
  some base of zero neighborhoods. We say that a net $(x_\gamma)$ is
  \term{$U$-bounded} if it is contained in $\alpha U$ for some $\alpha>0$,
  and $x_\gamma\xrightarrow{\su}0$ if for every $\alpha>0$ there exits
  $\gamma_0$ such that $x_\gamma\in\alpha U$ whenever
  $\gamma>\gamma_0$.

\begin{numbered}  \lab{n:nb-convergence}
  Suppose $T$ is nb-bounded, then $T(U)$ is bounded for some zero
  neighborhood $U$. Obviously $x_\gamma\xrightarrow{\su}0$ implies
  $Tx_\gamma\to 0$. The converse implication is also valid: if $T$
  maps $U$-convergent sequences to convergent sequences, then $T$ has
  to be nb-bounded and the set $T(U)$ is bounded. Indeed, if $T(U)$ is
  unbounded, then there is a zero neighborhood $V$ in $Y$ such that
  $V$ does not absorb $T(U)$. Then for every $n\ge 1$ there exists
  $y_n\in T(U)\setminus nV$. Suppose $y_n=Tx_n$ for some $x_n\in U$,
  then $\frac{x_n}{n}\xrightarrow{\su}0$, but
  $T(\frac{x_n}{n})=\frac{y_n}{n}\notin V$, so that $T(\frac{x_n}{n})$
  does not converge to zero.
\end{numbered}

\ssection{Normed, quasinormed, and seminormed spaces}

Next, we discuss bounded operators in some particular topologies.
Notice that every normed, seminormed, or quasinormed vector space is
locally bounded. Therefore bb-bound\-ed\-ness, continuity,
nn-bound\-ed\-ness and nb-bound\-ed\-ness coincide for operators on
such spaces.

\ssection{Locally convex topology}
%\marg{Would this work for quasi-seminorms?}

Similarly to the norm of an operator on a Banach space, we introduce
the seminorm of an operator on a seminormed space.

\begin{definition} \lab{d:pqseminorm}
  Let $T$ be an operator on a seminormed vector space $(X,p)$.
  As in the case with normed spaces,
  $p$ generates an \term{operator seminorm} $p(T)$ defined by 
  \[
    p(T) = \sup\limits_{p(x)\neq 0} \frac{p(Tx)}{p(x)}.
  \]
  More generally, let $S\colon X\to Y$ be a linear operator between
  two seminormed spaces $(X,p)$ and $(Y,q)$. Then we define a \term{mixed
  operator seminorm}\lab{a:msn} associated with $p$ and $q$ via
   \[
    \msn{p}{q}(S) = \sup\limits_{p(x)\neq 0} \frac{q(Sx)}{p(x)}.
   \]
\end{definition}

The seminorm $\msn{p}{q}(S)$ is a measure of how far in the seminorm $q$
the points of the $p$-unit ball can go under $S$. Notice, that
$p(T)$ and $\msn{p}{q}(S)$ may be infinite.
Clearly, if $T$ is an operator on a seminormed space $(X,p)$, then
$\msn{p}{p}(T)=p(T)$.

\begin{lemma} \lab{l:pqseminorm}
  If $S\colon X\to Y$ is an operator between two seminormed spaces
  $(X, p)$ and $(Y,q)$, then
  \begin{enumerate}
    \item $\msn{p}{q}(S) = \sup\limits_{p(x)= 1} q(Sx) =
          \sup\limits_{p(x)\le 1} q(Sx)$; \lab{i:os1}
    \item $q(Sx) \le \msn{p}{q}(S) p(x)$
            whenever $\msn{p}{q}(S)<\infty$. \lab{i:os2}
  \end{enumerate}
\end{lemma}

\begin{proof}
  The first equality in (\ref{i:os1}) follows immediately from the
  definition of $p(T)$. We obviously have
  \[
    \sup\limits_{p(x)=1} q(Sx) \le
    \sup\limits_{p(x)\le 1} q(Sx).
  \]
  In order to prove the opposite inequality, notice that if $0 < p(x) \le 1$, then
  $q(Sx) \le \frac{q(Sx)}{p(x)} \le \msn{p}{q}(S)$.
  Thus, it is left to show that $p(x) = 0$ implies
  $q(Sx) \le \msn{p}{q}(S)$. Pick any $z$ with $p(z) > 0$,
  then 
  \begin{displaymath}
    \textstyle
    p(\frac{z}{n})=p(x+\frac{z}{n}-x)\le
    p(x+\frac{z}{n})+p(x)=
    p(x+\frac{z}{n})\le p(x)+p(\frac{z}{n})=\frac{p(z)}{n},
  \end{displaymath}
  so that $p(x+\frac{z}{n})=p(\frac{z}{n})\in(0,1)$ for $n>p(z)$.
  Further, since $Sx+\frac{Sz}{n}$ converges to $Sx$ we have
  \[
    q(Sx)=\lim\limits_{n\to\infty}q\textstyle(Sx+\frac{Sz}{n})\displaystyle\le
    \lim\limits_{n\to\infty}
    \frac{q\left(S(x+\frac{z}{n})\right)}{p(x+\frac{z}{n})}
    \le \msn{p}{q}(S).
  \]
  
  Finally, (\ref{i:os2}) follows directly from the definition if
  $p(x)\neq 0$. In the case when $p(x)=0$, again pick any $z$ with
  $p(z)>0$, then $p(x+\frac{z}{n})\neq 0$ and
  \[
    q(Sx)=\lim\limits_{n\to\infty}q\textstyle(Sx+\frac{Sz}{n})\displaystyle
    =\lim\limits_{n\to\infty}q\left(\textstyle S(x+\frac{z}{n})\displaystyle\right)\le
    \lim\limits_{n\to\infty}\msn{p}{q}(S)p\textstyle(x+\frac{z}{n})
    =0.
  \]

\end{proof}

\begin{corollary}
  If $T$ is an operator on a seminormed space $(X, p)$, then
  \begin{enumerate}
    \item $p(T) = \sup\limits_{p(x)= 1} p(Tx)=\sup\limits_{p(x)\le 1} p(Tx)$;
    \item $p(Tx) \le p(T) p(x)$ whenever $p(T)<\infty$.
  \end{enumerate}
\end{corollary}

\bigskip

The following propositions characterize continuity and boundedness of
an operator on a locally convex space in terms of operator seminorms.
We assume that $X$ and $Y$ are locally convex spaces with generating
families of seminorms $\mathcal P$ and $\mathcal Q$ respectively.

\begin{proposition}
  Let $S$ be an operator from $X$ to $Y$, then $S$
  is continuous if and only if for every $q\in\mathcal Q$ there exists
  $p\in\mathcal P$ such that $\msn{p}{q}(S)$ is finite.
\end{proposition}

\begin{proposition} \lab{p:nn-smn}
  An operator $T$ on $X$ is nn-bounded if and only if $p(T)$ is finite
  for every $p\in\mathcal P$, or, equivalently, if $T$ maps
  $p$-bounded sets to $p$-bounded sets for every $p$ in some
  generating family $\mathcal P$ of seminorms.
\end{proposition}

\begin{proposition} \lab{p:nb-smn}
  Let $S\colon X\to Y$ be a linear operator, then the following are
  equivalent:
  \begin{enumerate}
    \item \lab{i:nb-smn1} $S$ is nb-bounded; 
    \item \lab{i:nb-smn2} $S$ maps $p$-bounded sets into bounded
      sets for some $p\in\mathcal P$;
    \item \lab{i:nb-smn3} There exists $p\in\mathcal P$ such that
             $\msn{p}{q}(S)<\infty$ for every $q\in\mathcal Q$. 
  \end{enumerate}
\end{proposition}

Since the balanced convex hall of a bounded set in a locally convex
space is again bounded, we also have the following characterization.

\begin{proposition} \lab{p:bb-smn}
  An operator $S\colon X\to Y$ is bb-bounded if and only if $\msn{p}{q}(S)$
  whenever $q\in\mathcal Q$ and $p$ is the Minkowski functional of a
  convex balanced bounded set.
\end{proposition}

\ssection{Operator topologies}

For each of the five classes of operators, we introduce an
appropriate natural operator topology. The class of all linear
operators between two topological vector spaces will be usually
equipped with the \term{strong} operator topology. Recall, that a
sequence $(S_n)$ of operators from $X$ to $Y$ is said to
\term{converge strongly} or \term{pointwise} to a map $S$ if $S_nx\to
Sx$ for every $x\in X$.  Clearly, $S$ will also be a linear operator.

\medskip

The class of all bb-bounded operators will usually be equipped with
the topology of \term{uniform convergence on bounded sets}. Recall,
that a sequence $(S_n)$ of operators is said to converge to zero
\term{uniformly} on $A$ if for each zero neighborhood $V$ in $Y$ there
exists an index $n_0$ such that $S_n(A)\subseteq V$ for all $n>n_0$.
We say that $(S_n)$ converges to $S$ uniformly on bounded sets if
$(S_n-S)$ converges to zero uniformly on bounded sets. Recall also
that a family $\mathcal G$ of operators is called \term{uniformly
  bounded} on a set $A\subseteq X$ if the set $\bigcup_{S\in{\mathcal
    G}}S(A)$ is bounded in $Y$.

\begin{lemma} \lab{l:bb-closed}
  If a sequence $(S_n)$ of bb-bounded operators converges uniformly on
  bounded sets to an operator $S$, then $S$ is also bb-bounded.
\end{lemma}
\marg{Exercise?}

\begin{proof}
  Fix a bounded set $A$. Since $S-S_n$ converges to zero uniformly on
  bounded sets then for every base zero neighborhood $V$ there exists
  an index $n_0$ such that $(S_n-S)(A)\subseteq V$ whenever $n\ge
  n_0$. This yields $S(A)\subseteq S_n(A)+V\subseteq \gamma V$ since
  $S_n(A)$ is bounded. Thus, $S(A)$ is bounded for every bounded set
  $A$, so that $S$ is bb-bounded.
\end{proof}

\medskip

The class of all continuous operators will be usually equipped with
the topology of \term{equicontinuous convergence}. Recall, that a
family $\mathcal G$ of operators from $X$ to $Y$ is called
\term{equicontinuous} if for each zero neighborhood $V$ in $Y$ there
is a zero neighborhood $U$ in $X$ such that $S(U)\subseteq V$ for
every $S\in\mathcal G$. We say that a sequence $(S_n)$ converges to zero
\term{equicontinuously} if for each zero neighborhood $V$ in $Y$ there
is a zero neighborhood $U$ in $X$ such that for every $\varepsilon>0$
there exists an index $n_0$ such that
$S_n(U)\subseteq \varepsilon V$ for all $n>n_0$.

\begin{lemma} \lab{l:cont-closed}
  If a sequence $S_n$ of continuous operators converges equicontinuously
  to $S$, then $S$ is also continuous.
\end{lemma}
\marg{Exercise?}

\begin{proof}
  Fix a zero neighborhood $V$, there exist zero neighborhoods $V_1$
  and $U$ and an index $n_0$ such that $V_1+V_1\subseteq V$ and
  $(S_n-S)(U)\subseteq V_1$ whenever $n>n_0$. Fix $n>n_0$. The
  continuity of $S_n$ guarantees that there exists a zero neighborhood
  $W\subseteq U$ such that $S_n(W)\subseteq V_1$. Since
  $(S_n-S)(W)\subseteq V_1$, we get $S(W)\subseteq S_n(W)+V_1\subseteq
  V_1+V_1\subseteq V$, which shows that $S$ is continuous.
\end{proof}

\medskip

The class of all nn-bounded operators will be usually equipped with
the topology of \term{nn-convergence}, defined as follows. We will
call a collection $\mathcal G$ of operators \term{uniformly
  nn-bounded} if there exists a base ${\mathcal N}_0$ of zero
neighborhoods such that for every $U\in{\mathcal N}_0$ there exists a
positive real $\beta$ such that $S(U)\subseteq\beta U$ for each
$S\in{\mathcal G}$. We say that a sequence $(S_n)$ \term{nn-converges}
to zero if there is a base ${\mathcal N}_0$ of zero neighborhoods such
that for every $U\in{\mathcal N}_0$ and every $\varepsilon>0$ we have
$S_n(U)\subseteq \varepsilon U$ for all sufficiently large $n$.

\begin{question}
  Is the class of all nn-bounded operators closed relative to
  nn-convergence? 
\end{question}

\medskip

Finally, the class of all nb-bounded operators will be usually
equipped with the topology of uniform convergence on a zero
neighborhood.

\begin{example}
  {\it The class of nb-bounded operators is not closed in the topology
    of uniform convergence on a zero neighborhood.} 
  Let $X={\mathbb R}^{\mathbb N}$, the space of all real sequences
  with the topology of coordinate-wise convergence.  Let $T_n$ be the
  projection on the first $n$ components. Clearly, every $T_n$ is
  nb-bounded because it maps the zero neighborhood
  $U_n=\bigl\{(x_i)_{i=1}^\infty\mid \abs{x_i}<1 \text{ for
    }i=1,\dots,n\bigr\}$ to a bounded set. On the other hand, $(T_n)$
  converges uniformly on $X$ to the identity operator, while the
  identity operator on $X$ is not nb-bounded.
\end{example}

\section{Spectra of an operator} \lab{s:spectra}

Recall that if $T$ is a continuous operator on a Banach space, then
its \term{resolvent set} $\rho(T)$ is the set of all
$\lambda\in\mathbb C$ such that the \term{resolvent operator}
$R_\lambda=(\lambda I-T)^{-1}$ exists, while the \term{spectrum} of
$T$ is defined by $\sigma(T)={\mathbb C}\setminus\rho(T)$. The Open
Mapping Theorem guarantees that if $R_\lambda$ exists then it is
automatically continuous. Now, if $T$ is an operator on an arbitrary
topological vector space and $\lambda\in\mathbb C$ then the algebraic
inverse $R_\lambda=(\lambda I-T)^{-1}$ may exist but not be
continuous, or may be continuous but not nb-bounded, etc. In order to
treat all these cases properly we introduce the following definitions.

\begin{definition} \lab{d:spectra}\lab{a:spectra}\lab{a:rhos}
  Let $T$ be a linear operator on a topological vector space. We
  denote the set of all scalars $\lambda\in\mathbb C$ for which
  $\lambda I-T$ is invertible in the algebra of linear operators by
  $\rho^l(T)$. We say that $\lambda\in\rho^{bb}(T)$ (respectively
  $\rho^c(T)$ or $\rho^{nn}(T)$) if the inverse of
  $\lambda I-T$ is bb-bounded (respectively continuous or nn-bounded).
  Finally, we say that $\lambda\in\rho^{nb}(T)$ if the inverse of
  $\lambda I-T$ belongs to the unitalization of the algebra of
  nb-bounded operators, i.e., when $(\lambda I-T)^{-1}=\alpha I+S$ for
  a scalar $\alpha$ and an nb-bounded operator $S$.

  The spectral sets $\sigma^l(T)$, $\sigma^{bb}(T)$,
  $\sigma^c(T)$, $\sigma^{nn}(T)$, and $\sigma^{nb}(T)$ are defined to
  be the complements of the resolvent sets $\rho^l(T)$,
  $\rho^{bb}(T)$,
  $\rho^c(T)$, $\rho^{nn}(T)$, and $\rho^{nb}(T)$ respectively.%
  \footnote{
    We use superscripts in order to avoid confusion with
    $\sigma_c(T)$, which is commonly used for continuous spectrum.}
  We will denote the (left and right) inverse of $\lambda I-T$
  whenever it exists by $R_\lambda$.
\end{definition}

\begin{numbered}  \lab{n:spectra-order}
  It follows immediately from Proposition~\ref{p:hierarchy} that
  $\sigma^l(T)\subseteq\sigma^{bb}(T)\subseteq\sigma^c(T)
  \subseteq\sigma^{nn}(T)\subseteq\sigma^{nb}(T)$. It follows from the
  Open Mapping Theorem that for a continuous operator $T$ on a Banach
  space all the introduced spectra coincide with the usual spectrum
  $\sigma(T)$. Since the Open Mapping Theorem is still valid on
  Fr\'echet spaces, we have $\sigma^l(T)=\sigma^{bb}(T)=\sigma^c(T)$
  for a continuous operator $T$ on a Fr\'echet space.
\end{numbered}

\begin{numbered} \lab{n:lb:spectra:equal}
  If $T$ is an operator on a locally bounded space $(X,U)$, then
  by~\ref{n:lb-equal} bb-boundedness of $T$ is equivalent to
  nb-boundedness, so that
  $\sigma^{bb}(T)=\sigma^c(T)=\sigma^{nn}(T)=\sigma^{nb}(T)$. We will
  denote this set by $\sigma_\su(T)$ to avoid ambiguity. Spectral
  theory of continuous operators on quasi-Banach spaces was developed
  in~\cite{Gramsch:66}.
\end{numbered}

\begin{numbered}  \lab{n:non-lb:sigma-nb:triv}
  There are several reasons why we define $\sigma^{nb}$ in a slightly
  different fashion than the other spectra. Namely, for $\lambda$ to
  be in $\rho^{nb}(T)$ we require $(\lambda I-T)^{-1}$ be not just
  nb-bounded, but be nb-bounded up to a multiple of the identity
  operator.  On one hand, this is the standard way to define the
  spectrum of an element in a non-unital algebra, and we know that the
  algebra of nb-bounded operators is unital only when the space is
  locally bounded. On the other hand, if we defined $\rho^{nb}(T)$ as
  the set of all $\lambda\in\mathbb C$ for which $(\lambda I-T)^{-1}$
  is nb-bounded, then we wouldn't have gotten any deep theory because
  $(\lambda I-T)^{-1}$ is almost never nb-bounded when the space is
  not locally bounded.
  
  Indeed, suppose that $X$ is not locally bounded, $T$ is a bb-bounded
  operator on $X$, and $\lambda\in\mathbb C$. Then $R_\lambda=(\lambda
  I-T)^{-1}$ cannot be nb-bounded, because in this case $I=(\lambda
  I-T)R_\lambda$ would be nb-bounded by~\ref{n:algebra} as a product
  of a bb-bounded and an nb-bounded operators. But we know that $I$ is
  not nb-bounded because $X$ is not locally bounded.
  
  We will see in Proposition~\ref{p:nb-not-invertible} that in a
  locally convex but non locally bounded space nb-bounded operators
  are never invertible, which implies that in such spaces $(\lambda
  I-T)^{-1}$ is not nb-bounded for any linear operator $T$.
\end{numbered}

\begin{numbered}  \lab{n:weak-spectra:equal}
  Next, let $T$ be a (norm) continuous operator on a Banach space,
  $\sigma(T)$ the usual spectrum of $T$, and let $\sigma^l(T)$,
  $\sigma^{bb}(T)$, $\sigma^c(T)$ be computed with respect to the weak
  topology. It is known that an operator on a Banach space is
  weak-to-weak continuous if and only if it is norm-to-norm
  continuous; therefore it follows that $\sigma^c(T)=\sigma(T)$.
  Furthermore, $\sigma^l(T)$ does not depend on the topology, so that
  it also coincides with $\sigma(T)$. Thus
  $\sigma^l(T)=\sigma^{bb}(T)=\sigma^c(T)=\sigma(T)$.
\end{numbered}

\section{Spectral radii of an operator} \lab{s:radii}

The spectral radius of a bounded linear operator $T$ on a Banach space
is usually defined via the Gelfand formula
$r(T)=\lim\limits_{n\to\infty}\sqrt[n]{\norm{T^n}}$. The formula
involves a norm and so makes no sense in a general topological vector
space. Fortunately, this formula can be rewritten without using a
norm, and then generalized to topological vector spaces. Similarly to
the situation with spectra, this generalization can be done in several
ways, so that we will obtain various types of spectral radii for an
operator on a topological vector space. We will show later that,
as with the Banach space case, there are some relations between
the spectral radii, the radii of the spectra, and the convergence of
the Neumann series of an operator on a locally convex topological
vector space. The content of this section may look technical at the
beginning, but later on the reader will see that all the facts lead to
a simple and natural classification.  We start with an almost obvious
numerical lemma.

\begin{lemma} \lab{l:numsr}
  If $(t_n)$ is a sequence in ${\mathbb R}^+\cup\{\infty\}$, then
  \[
    \limsup\limits_{n\to\infty}\sqrt[n]{t_n}=
    \inf\{\,\nu>0 \mid
      \lim\limits_{n\to\infty}\tfrac{t_n}{\nu^n}=0\,\}=
    \inf\{\,\nu>0 \mid
      \limsup\limits_{n\to\infty}\tfrac{t_n}{\nu^n}
      <\infty\,\}.
  \]
\end{lemma}

\begin{proof}
  Suppose $\limsup\limits_{n\to\infty}\sqrt[n]{t_n}=r$. If $0<\nu<r$,
  then $\sqrt[n_k]{t_{n_k}}>\mu>\nu$ for some $\mu$ and some
  subsequence $(t_{n_k})$, so that
  $\frac{t_{n_k}}{\nu^{n_k}}>\frac{\mu^{n_k}}{\nu^{n_k}}\to\infty$ as
  $k$ goes to infinity. It follows that
  $\limsup\lim_{n\to\infty}\frac{t_n}{\nu^n}=\infty$. On the other
  hand, if $r$ is finite and $\nu>r$ then $\sqrt[n]{t_n}<\mu<\nu$ for
  some $\mu$ and for all sufficiently large $n$. Then
  $\lim\limits_{n\to\infty}\frac{t_n}{\nu^n}\le
  \lim\limits_{n\to\infty}\frac{\mu^n}{\nu^n}=0$.
\end{proof}
  
This lemma implies that the spectral radius $r(T)$ of a (norm)
continuous operator $T$ on a Banach space equals the infimum of all
positive real scalars $\nu$ such that the sequence
$\bigl(\frac{T^n}{\nu^n}\bigr)$ converges to zero (or just is bounded)
in operator norm topology. This can be considered as an alternative
definition of the spectral radius, and can be generalized to any
topological vector space. Since for each of the five considered
classes of operators on topological vector spaces we introduced
appropriate concepts of convergent and bounded sequences, we arrive to
the following definition.

\begin{definition}  \lab{d:spradii}\lab{a:radii}
  Given a linear operator $T$ on a topological vector space $X$,
  define the following numbers:
  \begin{eqnarray*}
     r_l(T)&=&\inf\{\,\nu>0\mid \mbox{the sequence }\bigl(\tnun\bigr)
       \mbox{ converges strongly to zero}\,\};\\
     r_{bb}(T)&=&\inf\{\,\nu>0\mid\tnun\to 0
       \mbox{ uniformly on every bounded set}\,\};\\
     r_c(T)&=&\inf\{\,\nu>0\mid \tnun\to 0
       \mbox{ equicontinuously }\};\\
     r_{nn}(T)&=&\inf\{\,\nu>0\mid\bigl(\tnun\bigr)
       \mbox{ nn-converges to zero}\};\\
     r_{nb}(T)&=&\inf\{\,\nu>0\mid\tnun\to 0
       \mbox{ uniformly on some 0-neighborhood}\,\}.
   \end{eqnarray*}
\end{definition}

The following proposition explains the relations between the introduced
radii.

\begin{proposition} \lab{p:radii-order}
  If $T$ is a linear operator on a topological vector
  space $X$, then $r_l(T)\le r_{bb}(T)\le r_c(T)\le r_{nn}(T)\le r_{nb}(T)$.
\end{proposition}

\begin{proof}
  Let $T$ be a linear operator on a topological vector space $X$.
  Since every singleton is bounded then $r_l(T)\le r_{bb}(T)$. Next,
  assume $\nu>r_c(T)$, fix $\mu$ such that $r_c(T)<\mu<\nu$, then the
  sequence $(\frac{T^n}{\mu^n})$ converges to zero equicontinuously.
  Take a bounded set $A$ and a zero neighborhood $U$. There exists a
  zero neighborhood $V$ and a positive integer $N$ such that
  $\frac{T^n}{\mu^n}(V)\subseteq U$ whenever $n\ge N$.
  Also, $A\subseteq\alpha V$ for some $\alpha>0$, so that
  $\tnun(A)\subseteq\frac{\mu^n}{\nu^n}\frac{T^n}{\mu^n}(\alpha
  V)\subseteq\frac{\mu^n\alpha}{\nu^n}U\subseteq U$ for all
  sufficiently large $n$. It follows that the sequence $(\tnun)$
  converges to zero uniformly on $A$ and, therefore, $\nu\ge
  r_{bb}(T)$. Thus, $r_{bb}(T)\le r_c(T)$.
  
  To prove the inequality $r_c(T)\le r_{nn}(T)$ we let
  $\nu>r_{nn}(T)$. Then for some base ${\mathcal N}_0$ of zero
  neighborhoods and for every $V\in{\mathcal N}_0$ and $\varepsilon>0$
  there exists a positive integer $N$ such that
  $\tnun(V)\subseteq\varepsilon V$ for every $n\ge N$. Given a zero
  neighborhood $U$, we can find $V\in{\mathcal N}_0$ such that
  $V\subseteq U$. Then $\tnun(V)\subseteq\varepsilon
  V\subseteq\varepsilon U$ for every $n\ge N$, so that the sequence
  $(\tnun)$ converges to zero equicontinuously, and, therefore,
  $\nu\ge r_c(T)$.
  
  Finally, we must show that $r_{nn}(T)\le r_{nb}(T)$. Suppose that
  $\nu>r_{nb}(T)$, we claim that $\nu\ge r_{nn}(T)$. Take $\mu$ so
  that $\nu>\mu>r_{nb}(T)$. One can find a zero neighborhood $U$
  such that for every zero neighborhood $V$ there is a positive
  integer $N$ such that $\frac{T^n}{\mu^n}(U)\subseteq V$ for every
  $n\ge N$. Fix a base ${\mathcal N}_0$ of zero neighborhoods, and
  define a new base $\widetilde{\mathcal N}_0$ of zero neighborhoods
  via $\widetilde{\mathcal N}_0=\{mU\cap W\mid m\in{\mathbb
    N},\,W\in{\mathcal N}_0\}$. Let $V\in\widetilde{\mathcal N}_0$ and
  $\varepsilon>0$. Then $V=mU\cap W$ for some positive integer $m$ and
  $W\in{\mathcal N}_0$.  Then $\frac{T^n}{\mu^n}(V)\subseteq
  m\frac{T^n}{\mu^n}(U)\subseteq mV$ and for every sufficiently large
  $n$, so that $\frac{T^n}{\nu^n}(V)\subseteq
  \frac{\mu^n}{\nu^n}mV\subseteq\varepsilon V$, for each sufficiently
  large $n$, which implies $\nu\ge r_{nn}(T)$.
\end{proof}
  
The following lemma shows that, similarly to the case of Banach
spaces, one can use boundedness instead of convergence when defining
the spectral radii of an operator on a topological vector space. This
gives alternative ways of computing the radii, which are often more
convenient.

\begin{lemma}  \lab{l:radii-alt}
  Let $T$ be a linear operator on a topological vector space, then
  \begin{enumerate}
    \item \lab{li:rl-alt}
      $r_l(T)=\inf\bigl\{\,\nu>0\mid\bigl(\frac{T^nx}{\nu^n}\bigr)
      \text{\rm ~is bounded for every }x\,\bigr\};$
    \item \lab{li:rbb-alt}
      if $T$ is bb-bounded then\hfill\break
      $r_{bb}(T)=\inf\bigl\{\,\nu>0\mid\bigl(\tnun\bigr)
        \text{\rm ~is uniformly bounded on every bounded set}\,\bigr\};$
    \item  \lab{li:rc-alt}
      if $T$ is continuous then\hfill\break
      $r_c(T)=\inf\bigl\{\,\nu>0\mid\bigl(\tnun\bigr)
         \text{\rm ~is equicontinuous}\,\bigr\};$
    \item \lab{li:rnn-alt}
      if $T$ is nn-bounded then\hfill\break
      $r_{nn}(T)=\inf\bigl\{\,\nu>0\mid\bigl(\tnun\bigr)
         \text{\rm ~is uniformly nn-bounded}\,\bigr\};$
    \item \lab{li:rnb-alt}
      if $T$ is nb-bounded then\hfill\break
      $r_{nb}(T)=\inf\bigl\{\,\nu>0\mid\left(\tnun\right)\text{\rm ~is
         uniformly bounded on some 0-neighborhood}\,\bigr\}.$
  \end{enumerate}
  Moreover, in each of these cases it suffices to consider any tail of
  the sequence $\bigl(\tnun\bigr)$ instead of the whole sequence.
\end{lemma}

\begin{proof}
  To prove~(\ref{li:rl-alt}) let 
  \[
    r_l'(T)=\inf\bigl\{\,\nu>0\mid\bigl(\tnun x\bigr)
    \mbox{ is bounded for every }x\,\bigr\}.
  \]
  Since every convergent sequence is bounded, we certainly have
  $r_l(T)\ge r_l'(T)$. Conversely, suppose $\nu>r_l'(T)$, and take any
  positive scalar $\mu$ such that $\nu>\mu>r_l'(T)$. Then
  for every $x\in X$ the sequence $\frac{T^n}{\mu^n}x$ is bounded, and
  it follows that the sequence
  $\frac{T^nx}{\nu^n}=\frac{\mu^n}{\nu^n}\frac{T^nx}{\mu^n}$ converges
  to zero, so that $\nu\ge r_l(T)$ and, therefore $r_l'(T)\ge r_l(T)$.

  To prove~(\ref{li:rbb-alt}), suppose $T$ is bb-bounded, and let
  \[
    r_{bb}'(T)=\inf\bigl\{\,\nu>0\mid\bigl(\tnun\bigr)\mbox{ is uniformly
    bounded on every bounded set}\,\bigr\}.
  \]
  We'll show that $r_{bb}'(T)=r_{bb}(T)$. If
  $\bigl(\frac{T^n}{\nu^n}\bigr)$ converges to zero uniformly on every
  bounded set, then for each bounded set $A$ and for each zero
  neighborhood $U$ there exists a positive integer $N$ such that
  $\frac{T^n}{\nu^n}(A)\subseteq U$ whenever $n\ge N$. Also, since $T$
  is bb-bounded, then for every $n<N$ we have
  $\frac{T^n}{\nu^n}(A)\subseteq \alpha_n U$ for some $\alpha_n>0$.
  Therefore, if $\alpha=\max\{\alpha_1,\dots,\alpha_{N-1},1\}$, then
  $\frac{T^n}{\nu^n}(A)\subseteq\alpha U$ for every $n$, so that the
  sequence $\frac{T^n}{\nu^n}$ is uniformly bounded on $A$. Thus
  $\nu\ge r_{bb}'(T)$, so that $r_{bb}'(T)\le r_{bb}(T)$.
 
  Now suppose $\nu>r_{bb}'(T)$. There exists $\mu$ such that
  $\nu>\mu>r_{bb}'(T)$. The set
  $\bigcup_{n=1}^\infty\frac{T^n}{\mu^n}(A)$ is bounded for every
  bounded set $A$, so that for every zero neighborhood $U$ there
  exists a scalar $\alpha$ such that
  $\frac{T^n}{\mu^n}(A)\subseteq\alpha U$ for every $n\in\mathbb N$.
  Then $\tnun(A)\subseteq\frac{\mu^n\alpha}{\nu^n}U\subseteq U$ for
  all sufficiently large $n$. This means that the sequence
  $\bigl(\frac{T^n}{\nu^n}\bigr)$ converges to zero uniformly on $A$,
  and it follows that $\nu\ge r_{bb}(T)$.
  
  Further, if $T$ is bb-bounded, then any finite initial segment
  $(\tnun)_{n=0}^N$ is always uniformly bounded on every bounded
  set, so that a tail $(\tnun)_{n=N}^\infty$ is uniformly bounded on
  every bounded set if and only if the whole sequence
  $(\tnun)_{n=0}^\infty$ is uniformly bounded on every bounded set.

  The statements~(\ref{li:rc-alt}), (\ref{li:rnn-alt}),
  and~(\ref{li:rnb-alt}) can be proved in a similar way.
\end{proof}

\begin{numbered} \lab{n:lb:radii:equal}
  {\bf Locally bounded spaces.} If $T$ is a linear operator on a
  locally bounded topological vector space $(X,U)$, then it follows
  directly from Definition~\ref{d:spradii} that
  $r_{bb}(T)=r_c(T)=r_{nn}(T)=r_{nb}(T)$, because the corresponding
  convergences are equivalent. In this case we would denote each of
  these radii by $r_\su(T)$.
\end{numbered}

\ssection{Spectral radii via seminorms}

  The following proposition provides formulas for computing spectral radii
  of an operator on a locally convex space in terms of seminorms.
\marg{The same for semi-quasinorms?}

\begin{proposition} \lab{p:sr:lcs}
  If $T$ is an operator on a locally convex space $X$ with a generating
  family of seminorms $\mathcal P$, then
  \begin{enumerate}
    \item $r_l(T)=\sup\limits_{p\in{\mathcal P},\ x\in X}
          \limsup\limits_{n\to\infty}\sqrt[n]{p(T^nx)}$;
    \item $r_{bb}(T)=\sup\limits_{p\in{\mathcal B},\,q\in{\mathcal P}}
          \limsup\limits_{n\to\infty}\sqrt[n]{\msn{p}{q}(T^n)}$, where 
          $\mathcal B$ is the collection of Minkowski functionals of
          all balanced convex bounded sets in $X$;
    \item $r_c(T)=\sup\limits_{q\in{\mathcal P}}\inf\limits_{p\in{\mathcal P}}
          \limsup\limits_{n\to\infty}\sqrt[n]{\msn{p}{q}(T^n)}$;
    \item $r_{nn}(T)=\inf\limits_{\mathcal Q}\sup\limits_{p\in{\mathcal Q}}
          \limsup\limits_{n\to\infty}\sqrt[n]{p(T^n)}$,
          where the infimum is taken over all generating families of seminorms;
    \item $r_{nb}(T)=\inf\limits_{p\in{\mathcal P}}
          \sup\limits_{q\in{\mathcal P}}
          \limsup\limits_{n\to\infty}\sqrt[n]{\msn{p}{q}(T^n)}$;
  \end{enumerate}
\end{proposition}

\begin{proof}
  It follows from the definition of $r_l(T)$ and Lemma~\ref{l:numsr} that
  \begin{multline*}
     r_l(T)=
     \inf\bigl\{\,\nu>0\mid \lim\limits_{n\to\infty}
     p\left(\tfrac{T^nx}{\nu^n}\right)=0
     \mbox{ for every }x\in X, p\in{\mathcal P}\,\bigr\}\\
     =\sup\limits_{x\in X,\ p\in\mathcal P}
     \inf\bigl\{\,\nu>0\mid \lim\limits_{n\to\infty}
     \tfrac{p(T^nx)}{\nu^n}=0\,\bigr\}
     =\sup\limits_{x\in X,\ p\in{\mathcal P}}
     \limsup\limits_{n\to\infty}\sqrt[n]{p(T^nx)}.
  \end{multline*}
  
  Similarly, since the balanced convex hull of a bounded set is bounded,
  \begin{multline*}
    r_{bb}(T)=\inf\bigl\{\,\nu>0\mid\forall\mbox{ bounded } A\quad
    \forall V\in{\mathcal N}_0\quad\exists N\in\mathbb N\quad
    \forall n\ge N\quad\tnun(A)\subseteq V\,\bigr\}\\
    =\inf\bigl\{\,\nu>0\mid\forall p\in{\mathcal B}\quad\forall
    q\in{\mathcal P}\quad\exists N\in\mathbb N\quad
    \forall n\ge N\quad\msn{p}{q}\bigl(\tnun\bigr)\le 1\,\bigr\}\\
    =\sup\limits_{p\in{\mathcal B},\,q\in{\mathcal P}}
    \inf\bigl\{\,\nu>0\mid\lim\limits_{n\to\infty}
    \tfrac{\msn{p}{q}(T^n)}{\nu^n}\le 1\,\bigr\}
    =\sup\limits_{p\in{\mathcal B},\,q\in{\mathcal P}}
    \limsup\limits_{n\to\infty}\sqrt[n]{\msn{p}{q}(T^n)}.
  \end{multline*}

  Let $U_p=\{\,x\in X\mid p(x)<1\,\}$ for every $p\in\mathcal P$.
  Then, rephrasing the definition of $r_c(T)$ and applying
  Lemma~\ref{l:numsr}, we have
  \begin{multline*}
    r_c(T)=\inf\bigl\{\,\nu>0\mid\forall q\in{\mathcal P}\ \exists
    p\in{\mathcal P}\ \forall\varepsilon>0\ \exists N\in{\mathbb N}
    \ \forall n\ge N\ \tnun(U_p)\subseteq\varepsilon U_q\,\bigr\}\\
    =\sup\limits_{q\in{\mathcal P}}\inf\limits_{p\in{\mathcal P}}
    \inf\bigl\{\,\nu>0\mid\forall\varepsilon>0\ \exists N\in{\mathbb N}
    \ \forall n\ge N\ \msn{p}{q}\bigl(\tnun\bigr)<\varepsilon\,\bigr\}\\
%    =\inf\{\,\nu>0\mid\forall q\in{\mathcal P}\ \exists
%    p\in{\mathcal P}\ \exists N\in{\mathbb N}\ \forall n\ge N\ 
%    \msn{p}{q}(T^n)\le\nu^n\,\}\\
    =\sup\limits_{q\in{\mathcal P}}\inf\limits_{p\in{\mathcal P}}
    \inf\bigl\{\,\nu>0\mid\lim\limits_{n\to\infty}
    \tfrac{\msn{p}{q}(T^n)}{\nu^n}=0\,\bigr\}
    =\sup\limits_{q\in{\mathcal P}}\inf\limits_{p\in{\mathcal P}}
    \limsup\limits_{n\to\infty}\textstyle\sqrt[n]{\msn{p}{q}(T^n)}.
  \end{multline*}

  Similarly,
  \begin{multline*}
    r_{nn}(T)
      =\inf\bigl\{\,\nu>0\mid\exists{\mathcal Q}\ \forall p\in{\mathcal Q}\ 
      \forall\varepsilon>0\ \exists N\in{\mathbb N}\ \forall n>N\ 
      \tnun(U_p)\subseteq\varepsilon U_p\,\bigr\}\\
      =\inf\limits_{\mathcal Q}\sup\limits_{p\in{\mathcal Q}}
      \inf\bigl\{\,\nu>0\mid\lim\limits_{n\to\infty}
      \tfrac{p(T^n)}{\nu^n}=0\,\bigr\}
      =\inf\limits_{\mathcal Q}\sup\limits_{p\in{\mathcal Q}}
      \limsup\limits_{n\to\infty}\sqrt[n]{p(T^n)}.
  \end{multline*}

  Finally,
  \begin{multline*}
    r_{nb}(T)
    =\inf\bigl\{\,\nu>0\mid\exists p\in{\mathcal P}\ \forall
    q\in{\mathcal P}\ \exists N\in{\mathbb N}\ \forall n>N\ 
    \tnun(U_p)\subseteq U_q\,\bigr\}\\
    =\inf\limits_{p\in{\mathcal P}}\sup\limits_{q\in{\mathcal P}}
    \inf\bigl\{\,\nu>0\mid\limsup\limits_{n\to\infty}
    \tfrac{\msn{p}{q}(T^n)}{\nu^n}\le 1\,\bigr\}
    =\inf\limits_{p\in\mathcal P}\sup\limits_{q\in{\mathcal P}}
    \limsup\limits_{n\to\infty}\textstyle\sqrt[n]{\msn{p}{q}(T^n)}.
  \end{multline*}
\end{proof}

\ssection{Some special properties of $r_c(T)$}

Continuity of an operator can be characterized in terms of
neighborhoods (the preimage of every neighborhood contains a
neighborhood) or, alternatively, in terms of convergence (every
convergent net is mapped to a convergent net). Analogously, though
defined in terms of neighborhoods, $r_c(T)$ can also be characterized
in terms of convergent nets. This approach was used by F.~Garibay and
R.~Vera in a series of papers~\cite{Garibay:97,Garibay:98,Vera:97}.
Recall that a net $(x_\alpha)$ in a topological vector space is said
to be \term{ultimately bounded} if every zero neighborhood absorbs
some tail of the net, i.e., for every zero neighborhood $V$ one can
find an index $\alpha_0$ and a positive real $\delta>0$ such that
$x_\alpha\in\delta V$ whenever $\alpha>\alpha_0$. As far as we know,
ultimately bounded sequences were first studied in~\cite{DeVito:71}
for certain locally-convex topologies. The relationship between
ultimately bounded nets and convergence of sequences of operators on
locally convex spaces was studied
in~\cite{Garibay:97,Garibay:98,Vera:97}. The following proposition
(which is, in fact, a version of~\cite[Corollary 2.14]{Vera:97}) shows
how $r_c(T)$ can be characterized in terms of the action of powers of
$T$ on ultimately bounded sequences. It also implies that $r_c(T)$
coincides with the number $\gamma(T)$ which was introduced
in~\cite{Garibay:97,Garibay:98,Vera:97} for a continuous operator on
locally convex spaces.

%\begin{lemma}  \lab{l:ubub-ec}
%  For a sequence of linear operators $(S_n)$ on a topological vector
%  space the following are equivalent:
%  \begin{enumerate}
%    \item The net $(S_nx_\alpha)_{n,\alpha}$ is ultimately bounded
%      whenever the net $(x_\alpha)$ is ultimately bounded;
%    \item The sequence $(S_n)$ converges to zero equicontinuously.
%  \end{enumerate}
%\end{lemma}

\begin{proposition} \lab{p:rc-convergence}
  Let $T$ be a linear operator on a topological vector space $X$, then
  \begin{multline*}
    r_c(T)=
    \inf\bigl\{\nu>0\mid\lim\limits_{n,\alpha}\tfrac{T^n}{\nu^n}x_\alpha=0
      \text{ whenever $(x_\alpha)$ is ultimately bounded}\,\bigr\}\\
     =\inf\bigl\{\nu>0\mid\bigl(\tfrac{T^n}{\nu^n}x_\alpha\bigr)_{n,\alpha}
      \text{ is ultimately bounded whenever $(x_\alpha)$ is 
      ultimately bounded}\,\bigr\}.
  \end{multline*}
\end{proposition}

\begin{proof}
  To prove the first equality it suffices to show that $r_c(T)<1$ if
  and only if $\lim\limits_{n,\alpha}T^nx_\alpha=0$ whenever
  $(x_\alpha)$ is an ultimately bounded net. Suppose that $r_c(T)<1$,
  and let $V$ be a zero neighborhood. One can find a zero neighborhood
  $U$ such that for every $\varepsilon>0$ there exists $n_0\in\mathbb
  N$ such that $T^n(U)\subseteq\varepsilon V$ for each $n>n_0$. Let
  $(x_\alpha)$ be an ultimately bounded net. There exists an index
  $\alpha_0$ and a number $\delta>0$ such that $x_\alpha\in\delta U$
  whenever $\alpha>\alpha_0$. Then for $\varepsilon=\delta^{-1}$ one
  can find $n_0$ such that $T^n(U)\subseteq\delta^{-1}V$ for each
  $n>n_0$, so that $T^nx_\alpha\in\delta T^n(U)\subseteq V$ whenever
  $\alpha>\alpha_0$ and $n>n_0$. This means that
  $\lim\limits_{n,\alpha}T^nx_\alpha=0$.

  Conversely, suppose that $\lim\limits_{n,\alpha}T^nx_\alpha=0$ for
  each ultimately bounded net
  $(x_\alpha)$, and assume that $T^n$ does not converge
  equicontinuously to zero. Then there exists a zero
  neighborhood $V$ such that for every zero neighborhood $U$ one can
  find $\varepsilon_\su$ such that for every $m\in\mathbb N$ there exists
  $n_{\su,m}>m$ with $T^{n_{\su,m}}(U)\nsubseteq\varepsilon_\su V$.
  Then there exists $x_{\su,m}\in U$ such that
  \begin{equation}
    \label{eq:no-lim}
    T^{n_{\su,m}}x_{\su,m}\notin\varepsilon_\su V.    
  \end{equation}
  The collection of all zero neighborhood ordered by inclusion is a
  directed set, so that $(x_{\su,n})$ is an ultimately bounded net.
  Indeed, if $W$ is a zero neighborhood then 
  $x_{\su,n}\in W$ for each zero neighborhood $U\subseteq W$
  and every $n\in\mathbb N$. But it follows from~(\ref{eq:no-lim})
  that the net $\bigl(T^nx_{\su,m}\bigr)_{n,m,\su}$ does not converge
  to zero.

  To prove the second equality, let
  \begin{eqnarray*}
    \gamma_1&=&
    \inf\bigl\{\nu>0\mid\lim\limits_{n,\alpha}\tfrac{T^n}{\nu^n}x_\alpha=0
      \text{ whenever $(x_\alpha)$ is ultimately bounded}\,\bigr\}
    \mbox{ and}\\
     \gamma_2&=&
     \inf\bigl\{\nu>0\mid\bigl(\tfrac{T^n}{\nu^n}x_\alpha\bigr)_{n,\alpha}
      \text{ is ultimately bounded if $(x_\alpha)$ is 
      ultimately bounded}\,\bigr\}.
  \end{eqnarray*}
  Since every net which converges to zero is necessarily ultimately
  bounded, it follows that $\gamma_1\ge\gamma_2$. Now let
  $\nu>\gamma_2$, and let $(x_\alpha)$ be an ultimately bounded
  sequence. Suppose that $\gamma_2<\mu<\nu$, then
  $\bigl(\tfrac{T^n}{\mu^n}x_\alpha\bigr)_{n,\alpha}$ is ultimately
  bounded, that is, for each zero neighborhood $V$ there exists an
  indices $\alpha_0$ and $n_0$ and a positive $\varepsilon$ such that 
  $\tfrac{T^n}{\mu^n}x_\alpha\in\varepsilon V$ whenever
  $\alpha>\alpha_0$ and $n>n_0$. It follows that
  $\tfrac{T^n}{\nu^n}x_\alpha\in\frac{\mu^n\varepsilon}{\nu^n}V\subseteq V$
  for $\alpha>\alpha_0$ and all sufficiently large $n$. This implies
  that $\lim\limits_{n,\alpha}\tfrac{T^n}{\nu^n}x_\alpha=0$ so that
  $\nu\ge\gamma_1$. 
\end{proof}

%For a linear operator $T$ on a topological vector space $X$ we,
%following~\cite{Garibay:97,Garibay:98,Vera:97}, define the number
%$\gamma(T)$ as the infimum of the absolute value of those
%$\lambda\in\mathbb C$ for which
%$\lim\limits_{n,\alpha}\frac{T^nx_\alpha}{\lambda^n}=0$ whenever
%$(x_\alpha)$ is an ultimately bounded net. Since we assume that the
%base zero neighborhoods are balanced, we can, without loss of
%generality, consider only positive real $\lambda$'s.

%\begin{proposition}[{c.f.~\cite[Corollary~3.12]{Vera:97}}]  
%  Suppose that $T$ is a linear operator on a topological vector space
%  $X$, then $\gamma(T)=r_c(T)$.
%\end{proposition}

%\begin{proof}
%  If $\nu>\gamma(T)$ then
%  $\lim\limits_{n,\alpha}\frac{T^nx_\alpha}{\nu^n}=0$ for every
%  ultimately bounded net $(x_\alpha)$, so that by Lemma~\ref{l:ubub-ec}
%  the sequence $\tnun$ converges to zero equicontinuously. It follows
%  that $\nu\ge r_c(T)$, and thus $\gamma(T)\ge r_c(T)$.

%  Conversely, if $\nu>r_c(T)$, take a positive scalar $\mu$ such that
%  $\nu>\mu>r_c(T)$. Then $\frac{T^n}{\mu^n}$ converges to zero
%  equicontinuously. Fix an ultimately bounded net $(x_\alpha)$, then by
%  Lemma~\ref{l:ubub-ec} the net
%  $(\frac{T^nx_\alpha}{\mu^n})_{n,\alpha}$ is ultimately bounded, so
%  that the net $(\frac{T^nx_\alpha}{\nu^n})_{n,\alpha}$ converges to
%  zero. Thus $\nu\ge\gamma(T)$, so that $r_c(T)\ge\gamma(T)$.
%\end{proof}

\begin{question}
  Are there similar ways for computing
  $r_l(T)$, $r_{bb}(T)$, $r_{nn}(T)$, and $r_{nb}(T)$ in terms of
  nets?
\end{question}

Proposition~\ref{p:rc-convergence} enables us to prove some important
properties of $r_c$. The following lemma is analogous to Lemma~3.13
of~\cite{Vera:97}.

\begin{lemma} \lab{l:rc-prod}
  If $S$ and $T$ are two commuting linear operators on a topological
  vector space $X$ such that $r_c(S)$ and $r_c(T)$ are finite, then
  $r_c(ST)\le r_c(S)r_c(T)$.
\end{lemma}

\begin{proof}
  Suppose $\mu>r_c(S)$ and $\nu>r_c(T)$ and let $(x_\alpha)$ be an
  ultimately bounded net in~$X$. Then the net
  $(\frac{T^nx_\alpha}{\nu^n})_{n,\alpha}$ is ultimately bounded by
  Proposition~\ref{p:rc-convergence}. By applying
  Proposition~\ref{p:rc-convergence} again we conclude that
  $(\frac{S^mT^nx_\alpha}{\mu^m\nu^n})_{m,n,\alpha}$ converges to
  zero. In particular,
  $\lim\limits_{n,\alpha}\frac{(ST)^nx_\alpha}{(\mu\nu)^n}
  =\lim\limits_{n,\alpha}\frac{S^nT^nx_\alpha}{\mu^n\nu^n}=0$, and
  applying Proposition~\ref{p:rc-convergence} one more time we get
  $\mu\nu>r_c(ST)$.
\end{proof}

\begin{theorem} \lab{t:rc-sum}
  If $S$ and $T$ are two commuting continuous operators on a locally
  convex space $X$ then $r_c(S+T)\le r_c(S)+r_c(T)$.
\end{theorem}

\begin{proof}
  Assume without loss of generality that both $r_c(S)$ and $r_c(T)$
  are finite.  Suppose that $\eta>r_c(S)+r_c(T)$ and take
  $\mu>r_c(S)$ and $\nu>r_c(T)$ such that $\eta>\mu+\nu$. Let
  $(x_\alpha)$ be an ultimately bounded net in $X$. By
  Proposition~\ref{p:rc-convergence} it suffices to show that
  $\lim\limits_{n,\alpha}\frac{1}{\eta^n}(S+T)^nx_\alpha=0$. Notice
  that the net $\bigl(\frac{T^n}{\nu^n}x_\alpha\bigr)_{n,\alpha}$ is
  ultimately bounded. This implies that the net
  $\bigl(\frac{S^m}{\mu^m}\frac{T^n}{\nu^n}x_\alpha\bigr)_{m,n,\alpha}$
  converges to zero. Fix a seminorm $p$, then there exist indices
  $n_0$ and $\alpha_0$ such that $p(S^mT^nx_\alpha)<\mu^m\nu^n$
  whenever $m,n\ge n_0$ and $\alpha\ge\alpha_0$. Also, notice that we
  can split $\eta$ into a product of two terms $\eta=\eta_1\eta_2$
  such that $\eta_1>1$ while still $\eta_2>\mu+\nu$. Further, if
  $n>2n_0$ and $\alpha\ge\alpha_0$ then we have
  \begin{multline*}
    p\left(\tfrac{1}{\eta^n}(S+T)^nx_\alpha\right)\le\\
    \tfrac{1}{\eta^n}\sum\limits_{k=0}^{n_0}
      {\textstyle\binom{n}{k}}p\bigl(S^kT^{n-k}x_\alpha\bigr)+
    \tfrac{1}{\eta^n}\sum\limits_{k=n_0+1}^{n-n_0}
      {\textstyle\tbinom{n}{k}}p\bigl(S^kT^{n-k}x_\alpha\bigr)+
    \tfrac{1}{\eta^n}\!\!\!\sum\limits_{k=n-n_0+1}^n\!\!\!
      {\textstyle\tbinom{n}{k}}p\bigl(S^kT^{n-k}x_\alpha\bigr).
  \end{multline*}
  Since $\binom{n}{k}=\frac{(n-k+1)\cdots(n-1)\cdot n}{1\cdot
    2\cdots(k-1)\cdot k}\le n^k$ and
  $\sum_{k=0}^n\binom{n}{k}\mu^k\nu^{n-k}=(\mu+\nu)^n$, we have
  \begin{multline*}
    p\left(\tfrac{1}{\eta^n}(S+T)^nx_\alpha\right)\le\\
    \frac{n^{n_0}}{\eta^n}\sum\limits_{k=0}^{n_0}
      p\bigl(S^kT^{n-k}x_\alpha\bigr)+
    \frac{1}{\eta^n}\sum\limits_{k=n_0+1}^{n-n_0}
      {\textstyle\tbinom{n}{k}}\mu^k\nu^{n-k}+
    \frac{n^{n_0}}{\eta^n}\sum\limits_{k=n-n_0+1}^n
      p\bigl(S^kT^{n-k}x_\alpha\bigr)\\\le
    \frac{n^{n_0}}{\eta_1^n}\cdot\frac{1}{\eta_2^n}
      \sum\limits_{k=0}^{n_0}\Bigl(p\bigl(T^{n-k}S^kx_\alpha\bigr)+   
      p\bigl(S^{n-k}T^kx_\alpha\bigr)\Bigr)+
    \frac{(\mu+\nu)^n}{\eta^n}.
  \end{multline*}
  Notice that $\lim\limits_{n\to\infty}\frac{(\mu+\nu)^n}{\eta^n}=0$ and that
  $\lim\limits_{n\to\infty}\frac{n^{n_0}}{\eta_1^n}=0$.
  Since $T$ is continuous, the net $(T^kx_\alpha)_\alpha$ is
  ultimately bounded for every fixed $k$, so that
  $\lim\limits_{n,\alpha}\frac{1}{\eta_2^{n-k}}S^{n-k}T^kx_\alpha=0$. It
  follows that for every $k$ between $0$ and $n_0$ the expression
  $\frac{1}{\eta_2^n}p\bigl(S^{n-k}T^kx_\alpha\bigr)$ is uniformly bounded for
  all sufficiently large $n$ and $\alpha$. Similarly,
  for every $k$ between $0$ and $n_0$ the expression
  $\frac{1}{\eta_2^n}p\bigl(T^{n-k}S^kx_\alpha\bigr)$ is uniformly bounded for
  all sufficiently large $n$ and $\alpha$. Therefore there exist
  indices $n_1$ and $\alpha_1$ such that the finite sum
  \begin{displaymath}
    \frac{1}{\eta_2^n}\sum\limits_{k=0}^{n_0}
    \Bigl(p\bigl(T^{n-k}S^kx_\alpha\bigr)+p\bigl(S^{n-k}T^kx_\alpha\bigr)\Bigr) 
  \end{displaymath}
  is uniformly bounded for all $n\ge n_1$ and $\alpha\ge\alpha_1$.
  It follows that 
  $\lim\limits_{n,\alpha}p\Bigl(\frac{1}{\eta^n}(S+T)^nx_\alpha\Bigr)=0$,
  so that $\eta>r_c(S+T)$.
\end{proof}

\begin{corollary} \lab{c:poly-rc}
  If $T$ is a continuous operator on a locally convex space with
  finite $r_c(T)$ then $r_c\bigl(P(T)\bigr)$ is finite for every
  polynomial $P(z)$.
\end{corollary}

\begin{definition}
  We say that a sequence $(x_n)$ in a topological vector space is
  \term{fast null} if $\lim\limits_{n\to\infty}\alpha^n x_n=0$ for
  every positive real $\alpha$.
\end{definition}

\begin{lemma} \lab{l:fast-null}
  If $T$ is a linear operator on a topological vector space with
  $r_c(T)<\infty$ then $(T^nx_n)$ is fast null whenever $(x_n)$ is
  fast null.
\end{lemma}

\begin{proof}
  Suppose $(x_n)$ is a fast null sequence in a topological vector
  space and $r_c(T)<\infty$. Let $\nu>r_c(T)$, the sequence $\nu^n
  \alpha^n x_n$ converges to zero, hence is ultimately bounded, then
  by Proposition~\ref{p:rc-convergence} we have
  \begin{displaymath}
    \lim\limits_{n\to\infty}\alpha^nT^nx_n=
    \lim\limits_{n\to\infty}\tfrac{T^n}{\nu^n}\nu^n\alpha^nx_n=0.    
  \end{displaymath}
\end{proof}

\section{Spectra and spectral radii}  \lab{s:s-sr}

It is well known that for a continuous operator $T$ on a Banach space
its spectral radius $r(T)$ equals the geometrical radius of the
spectrum
$\bigabs{\sigma(T)}=\sup\{\abs{\lambda}\mid\lambda\in\sigma(T)\}$.
Further, whenever $|\lambda|>r(T)$, the resolvent operator
$R_\lambda=(\lambda I-T)^{-1}$ is given by the Neumann series
$\sum_{i=0}^\infty \frac{T^i}{\lambda^{i+1}}$.  We are going to show
in the next five theorems that the spectral radii that we have
introduced are upper bounds for the actual radii of the correspondent
spectra, and that when $\abs{\lambda}$ is greater than or equal to any
of these spectral radii, then the Neumann series converges in the
correspondent operator topology to the resolvent operator.

In the following Theorems~\ref{t:r-l}--\ref{t:r-nb} we assume that $T$
is a linear operator on a sequentially complete locally convex space,
$\lambda$ is a complex number, and $R_\lambda$ is the resolvent of $T$
at $\lambda$ in the sense of Definition~\ref{d:spectra}.

\begin{theorem} \lab{t:r-l}
  If $\abs{\lambda}>r_l(T)$ then the Neumann series converges
  pointwise to a linear operator $R_\lambda^0$\lab{a:r0}, and
  $R_\lambda^0(\lambda I-T)=I$.  Moreover, if $T$ is continuous, then
  $R_\lambda^0=R_\lambda$ and $\bigabs{\sigma^l(T)}\le r_l(T)$.
\end{theorem}

\begin{proof}
  For any $\lambda\in\mathbb C$ such that $\abs{\lambda}>r_l(T)$ one
  can find $z\in\mathbb C$ such that $0<|z|<1$ and $\lambda z>r_l(T)$.
  Consider a point $x\in X$ and a base zero neighborhood $U$. Since by
  the definition of $r_l(T)$ the sequence $\bigl(\frac{T^nx}{(\lambda
    z)^n}\bigr)$ converges to zero, there exist a positive integer
  $n_0$, such that $\frac{T^nx}{(\lambda z)^n}\in U$ whenever $n\ge
  n_0$. Therefore, $\frac{T^nx}{\lambda^n}\in
  z^nU\subseteq\abs{z}^nU$ because $U$ is balanced. Thus, if
  $n\ge m\ge n_0$, then
  $\sum_{i=n}^m\frac{T^ix}{\lambda^i}\in\sum_{i=n}^m\abs{z}^iU
  \subseteq\bigl(\sum_{i=n}^m\abs{z}^i\bigr)U$ because $U$ is convex.
  Since $\abs{z}<1$, we have $\sum_{i=n}^m\abs{z}^i<1$ for
  sufficiently large $m$ and $n$, and so $\sum_{i=n}^m
  \frac{T^ix}{\lambda^i}\in U$ because $U$ is balanced.  Therefore
  $R_{\lambda,n}x=\frac{1}{\lambda}\sum_{i=0}^n
  \frac{T^ix}{\lambda^i}$\lab{a:res-n} is a Cauchy sequence and hence
  it converges to some $R_\lambda^0x$ because $X$ is sequentially
  complete.
  
  Clearly, $R_\lambda^0$ is a linear operator. Notice that
  $R_{\lambda,n}(\lambda x-Tx)=x-\frac{T^{n+1}x}{\lambda^{n+1}}$ for
  every $x$. As $n$ goes to infinity, the left hand side of this
  identity converges to $R_\lambda^0(\lambda x-Tx)$, while the right
  hand side converges to $x$. Thus it follows that
  $R_\lambda^0(\lambda I-T)=I$.
  
  Finally, notice that $R_{\lambda,n}$ commutes with $T$ for every
  $n$. Therefore, if $T$ is continuous, then
  \[
    R_\lambda^0Tx=\lim\limits_{n\to\infty}R_{\lambda,n}Tx=
    \lim\limits_{n\to\infty}TR_{\lambda,n}x=
    T(\lim\limits_{n\to\infty}R_{\lambda,n}x)=
    TR_\lambda^0 x
  \]
  for every $x$. This implies that $(\lambda
  I-T)R_\lambda^0=R_\lambda^0(\lambda I-T)=I$, so that $R_\lambda^0$ is
  the (left and right) inverse of $\lambda I-T$. This means that
  $R_\lambda^0=R_\lambda$ and $\lambda\in\rho^l(T)$. Thus,
  $\bigabs{\sigma^l(T)}\le r_l(T)$.
\end{proof}

\begin{theorem} \lab{t:r-bb}
  If $T$ is bb-bounded and $\abs{\lambda}>r_{bb}(T)$, then the
  Neumann series converges uniformly on bounded sets, and its sum
  $R_\lambda^0$ is bb-bounded. Moreover, if $T$ is continuous, then
  $R_\lambda^0=R_\lambda$ and $\bigabs{\sigma^{bb}(T)}\le r_{bb}(T)$.
\end{theorem}

\begin{proof}
  Suppose that
  $\abs{\lambda}>r_{bb}(T)$, then the sum $R_\lambda^0$ of the Neumann series
  exists by Theorem~\ref{t:r-l}. As in the proof of
  Theorem~\ref{t:r-l} we denote the partial sums of the Neumann series
  by $R_{\lambda,n}$. Fix $z\in\mathbb C$ such that $0<\abs{z}<1$ and
  $\lambda z>r_{bb}(T)$, and consider a bounded set $A$ and a closed
  base zero neighborhood $U$. Since $\frac{T^n}{(\lambda z)^n}$
  converges to zero uniformly on $A$, there exits $n_0\in\mathbb N$
  such that $\frac{T^n}{\lambda^n z^n}(A)\subseteq U$ for all $n>n_0$.
  Also, since $\abs{z}<1$, we can assume without loss of generality
  that $\sum_{i=n_0}^\infty\abs{z}^i<\abs{\lambda}$.  Then
  \begin{displaymath}
    \tfrac{1}{\lambda}\sum_{i=n+1}^m\tfrac{T^ix}{\lambda^i}\in
    \tfrac{1}{\lambda}\Bigl(\sum_{i=n+1}^m\abs{z}^i\Bigr)U
    \subseteq U    
  \end{displaymath}
  whenever $x\in A$ and $m>n>n_0$. Since $U$ is closed, we have
  \begin{displaymath}
    R_\lambda^0 x-R_{\lambda,n}x=
    \lim\limits_{m\to\infty} \tfrac{1}{\lambda}
    \sum_{i=n+1}^m\tfrac{T^ix}{\lambda^i}\in U
  \end{displaymath}
  for each $x\in A$ and $n>n_0$, so that $(R_\lambda^0
  -R_{\lambda,n})(A)\subseteq U$ whenever $n>n_0$. This shows that
  $R_{\lambda,n}$ converges to $R_\lambda^0$ uniformly on bounded
  sets. By Lemma~\ref{l:bb-closed} this implies that $R_\lambda^0$ is
  bb-bounded.

%  To show that $R_\lambda^0$ is
%  bb-bounded, fix $z\in\mathbb C$ such that $0<\abs{z}<1$ and $\lambda
%  z>r_{bb}(T)$, and consider a bounded set $A$ and a closed base zero
%  neighborhood $U$. Since $\frac{T^n}{(\lambda z)^n}$ converges to
%  zero uniformly on $A$, then there exits $N\in\mathbb N$ such that
%  $\frac{T^n}{\lambda^n z^n}(A)\subseteq U$ whenever $n>N$. On the
%  other hand, the operator $\sum_{i=0}^N\frac{T^i}{\lambda^i}$ is
%  bb-bounded, so that
%  $\sum_{i=0}^N\frac{T^i}{\lambda^i}(A)\subseteq\gamma U$ for some
%  scalar $\gamma>0$.  If $x\in A$ then
%  \[
%    R_{\lambda,n}x=\frac{1}{\lambda}\sum\limits_{i=0}^n
%    \tfrac{T^ix}{\lambda^i}
%    \in\frac{1}{\lambda}\Bigl(\gamma+\sum\limits_{i=N+1}^n\abs{z}^i\Bigr)U
%  \]
%  whenever $n>N$.
%  Let $\xi_n=\frac{1}{\lambda}(\gamma+\sum_{i=N+1}^n\abs{z}^i)$, then
%  $\xi_n^{-1}R_{\lambda,n}x\in U$ for every $n>N$. Since
%  $\xi_n\to\xi=\frac{1}{\lambda}(\gamma+\sum_{i=N+1}^\infty\abs{z}^i)>0$
%  and since $U$ is closed, we have
%  $\xi^{-1}R_\lambda^0x\in U$, so that $R_\lambda^0x\in\xi U$.
%  Therefore $R_\lambda^0(A)\in\xi U$ whence $R_\lambda^0$
%  is bb-bounded.
  
  Further, if $T$ is continuous, then by Theorem~\ref{t:r-l} we
  have $R_\lambda=R_\lambda^0$, so that
  $\lambda\in\rho^{bb}(T)$, whence it follows that $\bigabs{\sigma^{bb}(T)}\le
  r_{bb}(T)$.
\end{proof}

The next theorem is similar to Theorem~2.18 of~\cite{Vera:97}.

\begin{theorem} \lab{t:r-c}
  If $T$ is a continuous and $\abs{\lambda}>r_c(T)$, then the Neumann
  series converges equicontinuously to $R_\lambda$, and $R_\lambda$ is
  continuous. In particular, $\bigabs{\sigma^c(T)}\le r_c(T)$ holds.
\end{theorem}

\begin{proof}
  Let $\abs{\lambda}>r_c(T)$. It follows from Theorem~\ref{t:r-l} that
  the Neumann series converges to $R_\lambda$. Again, we denote the
  partial sums of the Neumann series by $R_{\lambda,n}$.  Let
  $z\in\mathbb C$ be such that $0<|z|<1$ and $\lambda z>r_c(T)$.  For
  a fixed closed zero neighborhood $U$ there exists a zero
  neighborhood $V$ such that $\frac{T^n}{\lambda^n z^n}(V)\subseteq U$
  for every $n\ge 0$.  Let $\varepsilon>0$, then
  $\sum_{i=n_0}^\infty\abs{z}^i<\varepsilon\abs{\lambda}$ for some
  $n_0$. Then
  \begin{displaymath}
    \tfrac{1}{\lambda}\sum_{i=n+1}^m\tfrac{T^ix}{\lambda^i}\in
    \tfrac{1}{\lambda}\Bigl(\sum_{i=n+1}^m\abs{z}^i\Bigr)\varepsilon U
    \subseteq U    
  \end{displaymath}
  whenever $x\in V$ and $m>n>n_0$. Since $U$ is closed, we have
  \begin{displaymath}
    R_\lambda x-R_{\lambda,n}x=
    \lim\limits_{m\to\infty} \tfrac{1}{\lambda}
    \sum_{i=n+1}^m\tfrac{T^ix}{\lambda^i}\in\varepsilon U
  \end{displaymath}
  for each $x\in V$ and $n>n_0$, so that $(R_\lambda
  -R_{\lambda,n})(V)\subseteq\varepsilon U$ whenever $n>n_0$. This shows that
  $R_{\lambda,n}$ converges to $R_\lambda$ equicontinuously, and
  Lemma~\ref{l:cont-closed} yields that $R_\lambda$ is continuous.

%  In order to show that
%  $R_\lambda$ is continuous, fix some closed base zero neighborhood $U$ and
%  let $z\in\mathbb C$ be such that $0<|z|<1$ and $\lambda z>r_c(T)$.
%  Then there exists a zero neighborhood $V$ such that
%  $\frac{T^n}{\lambda^n z^n}(V)\subseteq U$ for every $n\ge 0$, and
%  for every $x\in V$ we have
%  \[
%    R_{\lambda,n}x=\frac{1}{\lambda}\sum\limits_{i=0}^n
%    \tfrac{T^ix}{\lambda^i}\in
%    \frac{1}{\lambda}\Bigl(\sum\limits_{i=0}^n\abs{z}^i\Bigr)U=
%    \frac{1-\abs{z}^n}{\lambda(1-\abs{z})}U
%  \]
%  for every $n\ge 0$. It follows that
%  $\frac{\lambda(1-\abs{z})}{1-\abs{z}^n}R_{\lambda,n}x\in U$, and then
%  $\lambda(1-\abs{z})R_\lambda x\in U$ because $U$ is closed.  Let
%  $W=\frac{1}{\lambda(1-\abs{z})}V$, then $R_\lambda(W)\subseteq U$, which
%  shows that $R_\lambda$ is continuous.
\end{proof}

\begin{theorem} \lab{t:r-nn}
  If $T$ is nn-bounded and $\abs{\lambda}>r_{nn}(T)$, then the Neumann
  series nn-converges to $R_\lambda$ and $R_\lambda$ is nn-bounded. In
  particular, $\bigabs{\sigma^{nn}(T)}\le r_{nn}(T)$ holds.
\end{theorem}

\begin{proof}
  Let $\abs{\lambda}>r_{nn}(T)$. By Theorem~\ref{t:r-l} the Neumann
  series $\sum_{i=0}^\infty\frac{T^i}{\lambda^{i+1}}$ converges to
  $R_\lambda$. Again, we denote the partial sums of the Neumann series
  by $R_{\lambda,n}$. Fix some $z$ such that $0<\abs{z}<1$ and
  $\lambda z>r_{nn}(T)$. There exists a base ${\mathcal N}_0$ of
  closed convex zero neighborhoods such that for every $U\in{\mathcal
    N}_0$ there is a scalar $\beta>0$ such that $\frac{T^n}{(\lambda
    z)^n}(U)\subseteq\beta U$ for all $n\ge 0$. Fix $U\in{\mathcal
    N}_0$, then for each $n\ge 0$ we have $\frac{T^n}{\lambda^n
    z^n}(U)\subseteq\beta U$ for some $\beta>0$, so that
  $\frac{T^nx}{\lambda^n}\in\abs{z}^n\beta U$ whenever $x\in U$. It
  follows that
  \[
    R_{\lambda,n}x=
    \tfrac{1}{\lambda}\sum_{i=0}^n\tfrac{T^ix}{\lambda^i}
    \in\tfrac{\beta}{\lambda}\Bigl(\sum\limits_{i=0}^n\abs{z}^i\Bigr)U.
  \]
  Then $R_\lambda x\in\frac{\beta}{\lambda(1-\abs{z})}U$, so that
  $R_\lambda(U)\subseteq\frac{\beta}{\lambda(1-\abs{z})}U$, which
  implies that $R_\lambda$ is nn-bounded, and, therefore,
  $\bigabs{\sigma^{nn}(T)}\le r_{nn}(T)$ holds.
  
  Fix $\varepsilon>0$. Then $\sum_{i=N}^\infty\abs{z}^i<\abs{\lambda}$
  for some $N$.  Then for every $U\in{\mathcal N}_0$ we have
  \begin{displaymath}
    \tfrac{1}{\lambda}\sum_{i=n+1}^m\tfrac{T^ix}{\lambda^i}\in
    \tfrac{1}{\lambda}\Bigl(\sum_{i=n+1}^m\abs{z}^i\Bigr)\varepsilon U
    \subseteq U    
  \end{displaymath}
  whenever $x\in U$ and $N<n<m$. Since $U$ is closed, we have
  \begin{displaymath}
    R_\lambda x-R_{\lambda,n}x=
    \lim\limits_{m\to\infty} \tfrac{1}{\lambda}
    \sum_{i=n+1}^m\tfrac{T^ix}{\lambda^i}\in\varepsilon U
  \end{displaymath}
  for each $x\in U$ and $n>N$, so that $(R_\lambda
  -R_{\lambda,n})(U)\subseteq\varepsilon U$ whenever $N<n$. This shows that
  $R_{\lambda,n}$ nn-converges to $R_\lambda$.
\end{proof}

\begin{theorem} \lab{t:r-nb}
  If $T$ is nb-bounded and $\abs{\lambda}>r_{nb}(T)$, then the Neumann
  series converges to $R_\lambda$ uniformly on a zero neighborhood.
  Further, $\bigabs{\sigma^{nb}(T)}\le r_{nb}(T)$ holds.
\end{theorem}

\begin{proof}
  Let $\abs{\lambda}>r_{nb}(T)$. By Theorem~\ref{t:r-l} the Neumann
  series $\sum_{i=0}^\infty\frac{T^i}{\lambda^{i+1}}$ converges to
  $R_\lambda$. Since $r_{bb}(T)\le r_{nb}(T)$ then $R_\lambda$ is
  bb-bounded by Theorem~\ref{t:r-bb}. But then
  $\sum_{i=0}^\infty\frac{T^i}{\lambda^{i+1}}=\frac{1}{\lambda}I+
  \frac{1}{\lambda}R_\lambda T$. Notice that $R_\lambda T$ is
  nb-bounded as a product of a bb-bounded and an nb-bounded operators
  (see~\ref{n:algebra}).
  
  Suppose that $\abs{\lambda}>r_{nb}(T)$. Fix $z\in\mathbb C$ such
  that $0<\abs{z}<1$ and $\lambda z>r_{nb}(T)$, then the sequence
  $\bigl(\frac{T^n}{\lambda^n z^n}\bigr)$ converges to zero
  uniformly on some base zero neighborhood~$U$. We will show that the
  Neumann series converges uniformly on $U$. As in the proof of
  Theorem~\ref{t:r-l}, we denote the partial sums of the Neumann series
  by $R_{\lambda,n}$. Fix a closed base zero neighborhood $V$. Since
  $\bigl(\frac{T^n}{\lambda^n z^n}\bigr)$ converges to zero uniformly
  on $U$, there exits $n_0\in\mathbb N$ such that
  $\frac{T^n}{\lambda^n z^n}(U)\subseteq V$ for all $n>n_0$. Also,
  since $\abs{z}<1$, we can assume without loss of generality that
  $\sum_{i=n_0}^\infty\abs{z}^i<\abs{\lambda}$.  Then
  \begin{displaymath}
    \tfrac{1}{\lambda}\sum_{i=n+1}^m\tfrac{T^ix}{\lambda^i}\in
    \tfrac{1}{\lambda}\Bigl(\sum_{i=n+1}^m\abs{z}^i\Bigr)V
    \subseteq V
  \end{displaymath}
  whenever $x\in A$ and $m>n>n_0$. Since $V$ is closed, we have
  \begin{displaymath}
    R_\lambda x-R_{\lambda,n}x=
    \lim\limits_{m\to\infty} \tfrac{1}{\lambda}
    \sum_{i=n+1}^m\tfrac{T^ix}{\lambda^i}\in V
  \end{displaymath}
  for each $x\in U$ and $n>n_0$, so that $(R_\lambda
  -R_{\lambda,n})(U)\subseteq V$ whenever $n>n_0$.
\end{proof}

\medskip

In rest of this section we present some remarks on
Theorems~\ref{t:r-l}--\ref{t:r-nb}. In particular, we discuss the
conditions of sequential completeness and the local convexity and
consider several examples and special cases.

\begin{numbered}
  It is easy to see that each spectral radius is exactly the radius of
  convergence of the Neumann series in the correspondent operator
  convergence. Indeed, in each of Theorems~\ref{t:r-l}--\ref{t:r-nb}
  the convergence of the Neumann series implies that the terms of the
  series tend to zero. It follows that $\abs{\lambda}$ is greater than
  or equal to the corresponding spectral radius.
\end{numbered}

Clearly, if $X$ is a Banach space, then the norm topology on $X$ and
the weak${}^*$ topology on $X^*$ are sequentially complete. The weak
topology of $X$ is sequentially complete if $X$ is reflexive. Also, it
is known that the weak topologies of $\ell_1$ and of $L_1[0,1]$ are
sequentially complete. Since all these topologies are locally convex,
Theorems~\ref{t:r-l}--\ref{t:r-nb} are applicable to each of them.

\begin{numbered} \lab{n:MCP}
  {\bf Monotone convergence property.} Notice that if $T$ is a
  positive operator on a locally convex-solid vector lattice
  (i.e., a locally convex space which is also a vector lattice
    such that $\abs{x}\le\abs{y}$ implies $p(x)\le p(y)$ for every
    generating seminorm $p$)
  then we can substitute the sequential
  completeness condition in Theorems~\ref{t:r-l}--\ref{t:r-nb}
  by a weaker condition called
  \term{sequential monotone completeness property}: a locally
  convex-solid vector lattice is said to satisfy the sequential
  monotone completeness property if every monotone Cauchy sequence
  converges in the topology of $X$. For details,
  see~\cite{Aliprantis:78}. Indeed, we used the sequential
  completeness at just one point --- we used it in the proof of
  Theorem~\ref{t:r-l} to claim that since $R_{\lambda,n}x=
  \frac{1}{\lambda}\sum_{i=0}^n \frac{T^ix}{\lambda^i}$ is a Cauchy
  sequence, then it converges to some $R_\lambda x$. But if $T$ is
  positive, then $R_{\lambda,n}x^+$ and $R_{\lambda,n}x^-$ are
  increasing sequences, and the sequential monotone completeness
  property ensures the convergence.
\end{numbered}

\begin{numbered}
  {\bf Pointwise convergence.}
  It can be easily verified that the space of continuous functions on
  $[0,1]$ with pointwise convergence topology is not sequentially
  complete, the sequence $x_n(t)=t^n$ is a counterexample. The same
  counterexample shows that this space 
  does not have the monotone convergence property either.
  
  Consider the sequence spaces $\ell_p$ for $0 < p\le
  \infty$, $c$, $c_0$, and $c_{00}$ (the space of eventually vanishing
  sequences). None of these spaces is sequentially complete in the
  topology of coordinate-wise convergence: take the following sequence for a
  counterexample:
   \[
     x_n(i)=
     \left\{ 
       \begin{array}{ll}
         i & \mbox{if $i<n$};\\
         0 & \mbox{otherwise}.\\
       \end{array}
     \right.
   \]
   The same example shows that these spaces do not have the monotone
   convergence property either. Therefore neither of
   Theorems~\ref{t:r-l}--\ref{t:r-nb} or \ref{n:MCP} can be applied.
\end{numbered}

\begin{example}
  {\it Theorems~\ref{t:r-l}--\ref{t:r-nb} fail without sequential
    completeness.}
  Consider the space $c_0$ with the topology of coordinate-wise
  convergence.  Let $T$ be the forward shift operator on $c_0$, that
  is, $Te_k=e_{k+1}$, where $e_k$ is the $k$-th unit vector of $c_0$.
  Let $V$ be any base zero neighborhood, we can assume without loss of
  generality that $V=\{x\in\nolinebreak c_0\mid\abs{x_{i_1}}<1,\dots,
  \abs{x_{i_k}}<1\}$ where $i_1<i_2<\dots
  <i_k$ are positive integers. If $x\in U$ then $T^nx$ has zero
  components $1$ through $n$, in particular for every positive $\nu$
  we have $\frac{T^nx}{\nu^n}\in V$ whenever $n>i_k$. Therefore
  $\bigl(\tnun\bigr)$ converges uniformly on $c_0$ for every $\nu>0$,
  so that $r_{nb}(T)=0$. It follows from
  Proposition~\ref{p:radii-order} that
  $r_l(T)=r_{bb}(T)=r_c(T)=r_{nn}(T)=0$. On the other hand,
  $\sum_{n=1}^\infty T^ne_1$ diverges in $c_0$. Since $T$ is obviously
  continuous, this shows that Theorems~\ref{t:r-l}--\ref{t:r-nb}, do
  not hold in $c_0$. Thus, sequential completeness condition is
  essential in the theorems.
\end{example}

\begin{numbered}  \lab{n:ban:radii:equal}
  {\bf Banach spaces.} If $T$ is a (norm) continuous operator on a
  Banach space, then it follows from~\ref{n:spectra-order}
  and~\ref{n:lb:radii:equal} that
  $\sigma^l(T)=\sigma^{bb}(T)=\sigma^c(T)=\sigma^{nn}(T)=
  \sigma^{nb}(T)=\sigma(T)$ and
  $r_{bb}(T)=r_c(T)=r_{nn}(T)=r_{nb}(T)=r(T)$, where $\sigma(T)$ and
  $r(T)$ are the usual spectrum and the spectral radius of $T$.
  Further, it follows from Lemma~\ref{p:radii-order} that $r_l(T)\le
  r(T)$. On the other hand, since $r(T)=\bigabs{\sigma(T)}$, then
  $r(T)\le r_l(T)$ by Theorem~\ref{t:r-l}, so that $r_l(T)=r(T)$.
\end{numbered}

%\begin{numbered}
%  {\bf Weak topologies.} Suppose $X$ is a Banach space, and let
%  $(f_n)$ be a weak${}^*$-Cauchy sequence in $X^*$. Then $f_n(x)$ is
%  Cauchy for every $x\in X$, and therefore $f_n(x)$ converges to some
%  $f(x)$. By Banach-Steinhouse $f\in X^*$, and then $f$ is the
%  weak${}^*$-limit of $(f_n)$, so that weak${}^*$ topology is
%  sequentially complete. In particular, every reflexive Banach space
%  is weakly sequentially complete.
  
%\marg{Remove this, just state that $\ell_1$ and $L[0,1]$ are
%  sequentially complete?}
%  By Schur's Lemma weak convergence coincides with the norm
%  convergence for sequences in $\ell_1$. It follows that $\ell_1$ is weakly
%  sequentially complete. Finally, every KB-space is weakly
%  sequentially complete (see \cite[Theorem X.4.9]{KA}). In particular,
%  $L_1[0,1]$ is weakly sequentially complete.
%\end{numbered}

\begin{numbered}
  The following argument is a counterpart
  to~\ref{n:weak-spectra:equal}.  Let $T$ be a (norm) continuous
  operator on a Banach space $X$ and $r(T)$ the usual spectral radius
  of $T$, while $r_l(T)$ and $r_{bb}(T)$ be computed with respect to
  the weak topology of $X$. We claim that if the weak topology of $X$ is
  sequentially complete, then $r_l(T)=r_{bb}(T)=r(T)$. Indeed,
  $r(T)\le r_l(T)$ by~\ref{n:weak-spectra:equal} and
  Theorem~\ref{t:r-l} because $\sigma(T)=\sigma^l(T)$. In view of
  Proposition~\ref{p:radii-order} it suffices to show that
  $r_{bb}(T)\le r(T)$. Let $\nu>r(T)$, and let $A$ be a weakly bounded
  subset of $X$. Then $A$ is norm bounded, so that the sequence
  $\tnun$ converges to zero uniformly on $A$ in the norm topology. In
  particular, the set $\bigcup_{n=0}^\infty\tnun(A)$ is norm bounded,
  hence weakly bounded, so that $\nu>r_{bb}(T)$.
\end{numbered}

\ssection{Quasinilpotence}

Recall that a norm continuous operator $T$ on a Banach space $X$ is
said to be \term{quasinilpotent} if $r(T)=0$ or, equivalently, if
$\sigma(T)=\{0\}$. Quasinilpotent operators on Banach spaces have some
nice properties, therefore in the framework of topological vector
spaces it is interesting to study operators having some of their
spectra trivial or some of their spectral radii being zero. Notice,
for example, that it follows from Proposition~\ref{p:sr:lcs} that if
$T$ is an operator on a locally convex topological vector space, then
$r_l(T)=0$ if and only if
$\lim\limits_{n\to\infty}\sqrt[n]{p(T^nx)}=0$ for every seminorm $p$
in a generating family of seminorms and for every $x\in X$. Further,
if the space is in addition sequentially complete, then for such an
operator we would have $\sigma^l(T)=\{0\}$ by Theorem~\ref{t:r-l}.

Recall also that a norm continuous operator $T$ on a Banach space $X$
is said to be \term{locally quasinilpotent at a point} $x\in X$ if
$\lim\limits_{n\to\infty}\sqrt[n]{\norm{T^nx}}=0$. Using
Lemma~\ref{l:numsr}, the concept of local quasinilpotence can be
naturally generalized to topological vector spaces: an operator $T$ on
a topological vector space $X$ is said to be \term{locally
  quasinilpotent at a point} $x\in X$ if
$\lim\limits_{n\to\infty}\frac{T^nx}{\nu^n}=0$ for every $\nu>0$.  It
follows immediately from the definition of $r_l(T)$ that $r_l(T)=0$ if
and only if $T$ is locally quasinilpotent at every $x\in X$. It is
known that a continuous operator on a Banach space is quasinilpotent
if and only if it is locally quasinilpotent at every point. We see now
that this is just a corollary of~\ref{n:ban:radii:equal}.  The
following example shows that a similar result for general topological
vector spaces is not valid, that is, $r_l(T)$ may be equal to zero
without the other radii be equal to zero.

\begin{example}
  {\it A continuous operator with $r_l(T)=0$ but
  $r_{bb}(T)=r_c(T)=r_{nn}(T)=r_{nb}(T)=\infty$.}
  Consider the space of all bounded real sequences
  $\ell_\infty=\bigl\{x=(x_1,x_2,\dots)\mid \sup|x_k|<\infty\bigr\}$
  with the topology of coordinate-wise convergence. This topology can be
  generated by the family of coordinate seminorms $\{p_m\}_{m=1}^\infty$
  where $p_m(x)=\abs{x_m}$. Let $e_k$ denote the $k$-th unit vector in
  $\ell_\infty$.
  
  Define an operator $T\colon\ell_\infty\to\ell_\infty$ via
  $Te_k=\frac{(k-1)^{k-1}}{k^k}e_{k-1}$ if $k>1$, and $Te_1=0$. Then
  $T^ne_k= \frac{(k-n)^{k-n}}{k^k}e_{k-n}$ if $n<k$ and zero
  otherwise. Clearly $T$ is continuous. In order to show that
  $r_l(T)=0$ fix a positive real number $\nu$ and $x\in\ell_\infty$,
  then
  \begin{displaymath}
    \bigabs{\left(\tfrac{T^nx}{\nu^n}\right)_m}=
    \bigabs{\tfrac{m^m}{(m+n)^{m+n}\nu^n}x_{n+m}}\le
    \sup\limits_{n}\tfrac{m^m}{(m+n)^{m+n}\nu^n}\cdot
    \sup\limits_{n}\abs{x_n}<\infty
  \end{displaymath}
  It follows from Lemma~\ref{l:radii-alt}(\ref{li:rl-alt}) that
  $r_l(T)=0$.
 
  Now we show that $r_{bb}(T)=\infty$ by presenting a bounded set $A$ 
  in $\ell_\infty$ such that the sequence
  $\bigl(\frac{T^n}{\nu^n}\bigr)$ in not uniformly bounded on $A$ for
  every positive $\nu$. Let
  \begin{displaymath}
    A=\bigl\{\,x\in\ell_\infty\mid
      x_n\le(2n)^{2n}\text{ for all }n\ge 0\,\bigr\}.
  \end{displaymath}
  Then $(2n)^{2n}e_n\in A$ for each $n>0$ and
  $\bigl(\frac{T^{n-1}}{\nu^{n-1}}(2n)^{2n}e_n\bigr)_1=
  \frac{(2n)^{2n}}{n^n\nu^n}$ is unbounded. Then by
  Lemma~\ref{l:radii-alt}(\ref{li:rbb-alt}) we have
  $r_{bb}(T)=\infty$, and it follows from
  Proposition~\ref{p:radii-order} that
  $r_c(T)=r_{nn}(T)=r_{nb}(T)=\infty$.

\marg{Do this in details? Would take a couple paragraphs.}
  It is not difficult to show that $\sigma^l(T)=\{0\}$, while
  $\sigma^c(T)=\sigma^{nn}(T)=\sigma^{nb}(T)=\mathbb C$.
\end{example}

\ssection{Non-locally convex spaces}

We proved the key Theorems~\ref{t:r-l}--\ref{t:r-nb} for locally convex
spaces, but they are still valid for locally pseudo-convex spaces. The
local convexity of $X$ was used only once in the proof of
Theorem~\ref{t:r-l}, while Theorems~\ref{t:r-bb}--\ref{t:r-nb} used
Theorem~\ref{t:r-l}. Hence it would suffice to modify the proof of
Theorem~\ref{t:r-l} in such a way that it would work for locally
pseudo-convex spaces instead of locally convex.  Local convexity was
used in the proof of Theorem~\ref{t:r-l} to show that if
$\frac{T^nx}{(\lambda z)^n}\in U$ for all $n>n_0$ and some
$n_0\in\mathbb N$, then there exists $m_0\in\mathbb N$ such that
$\sum_{i=n}^m\frac{T^ix}{\lambda^i}\in U$ for all $m,n>m_0$.  (Recall
that $T$ is a linear operator, $\lambda,z\in\mathbb C$ such that
$0<\abs{z}<1$ and $\lambda z>r_l(T)$, $x\in X$, and $U$ is a base zero
neighborhood in $X$.)  If $X$ is locally pseudo-convex, then we can
assume that $U+U\subseteq\alpha U$ for some $\alpha>0$, so that
$(X,U)$ is a locally bounded space. Let $\norm{\cdot}$ be the
Minkowski functional of $U$, then (see~\cite[pages 3 and
6]{Kalton:84}) for any $x_1,x_2,\dots,x_n$ in $X$ we have
\begin{displaymath}
  \norm{x_1+...+x_k}\le 4^{\frac{1}{p}}\bigl(\norm{x_1}^p+\dots
    +\norm{x_k}^p\bigr)^{\frac{1}{p}},
\end{displaymath}
where $2^{\frac{1}{p}}=\alpha$. Notice that
$\norm{\frac{T^nx}{\lambda^n}}\le\abs{z}^n$ for all $n>n_0$.
Since $\abs{z}<1$, then there exists $m_0$
such that
$\sum_{i=n}^m\abs{z}^{ip}<\frac{1}{4}$ whenever $n,m>m_0$. But then
\begin{displaymath}
  \Bignorm{\sum\limits_{i=n}^m\tfrac{T^ix}{\lambda^i}}\le
  4^{\frac{1}{p}}\Bigl(\sum\limits_{i=n}^m\norm{
      \tfrac{T^ix}{\lambda^i}}^p\Bigr)^{\frac{1}{p}}\le
  4^{\frac{1}{p}}\Bigl(\sum\limits_{i=n}^m\abs{z}^{ip}\Bigr)^{\frac{1}{p}}<1,
\end{displaymath}
so that $\sum\limits_{i=n}^m\frac{T^ix}{\lambda^i}\in U$.

The following example shows that Theorems~\ref{t:r-l}, \ref{t:r-bb},
and~\ref{t:r-c} fail if we assume no convexity conditions at all.

\begin{example}
  {\it An operator on a complete non locally pseudo-convex space,
    whose spectral radii are 1, and whose Neumann series
    nevertheless diverges at $\lambda=2$.}  Let $X$ be the space of
  all measurable functions on $[0,1]$ with the topology of convergence
  in measure (which is not pseudo-convex). We identify the endpoints
  $0$ and $1$ and consider the interval as a circle. Fix an irrational
  $\alpha$ and define a linear operator $T$ on $X$ as the translation
  by $\alpha$, i.e., $(Tf)(t)=f(t-\alpha)$.  It is easy to see that
  $\frac{T^nf}{\nu^n}$ converges in measure to zero for every $f\in X$
  if and only if $\nu>1$. We conclude, therefore, that $r_l(T)=1$.
  Moreover, since the sets of the form
  $W_{\varepsilon,\delta}=\bigl\{f\in X
  \mid\mu(f>\varepsilon)<\delta\bigr\}$ form a zero neighborhood base
  for the topology of convergence in measure, and
  $T(W_{\varepsilon,\delta})\subseteq W_{\varepsilon,\delta}$, it
  follows that $r_{nn}(T)\le 1$. Then by
  Proposition~\ref{p:radii-order} we have
  $r_l(T)=r_{bb}(T)=r_c(T)=r_{nn}(T)=1$. Nevertheless, we are going to
  present a function $h\in X$ such that the Neumann series
  $\sum_{n=0}^\infty\frac{T^nh}{2^n}$ does not converge in measure,
  which means that the conclusions of
  Theorems~\ref{t:r-l}--\ref{t:r-nb} do not hold for this space

  For each $n=1,2,3,\dots$ one can find a positive integer $M_n$ such
  that the intervals $[k\alpha,k\alpha+\frac{1}{n}]$ (mod 1) for
  $k=1,\dots,M_n$ cover the circle. Let $s_n=\sum_{i=1}^nM_i$, and let
  $h$ be the step function taking value $2^{s_n}$ on the interval
  $(\frac{1}{n+1},\frac{1}{n}]$. If $s_{n-1}<k\le s_n$ for some
  positive integers $n$ and $k$, then on $[0,\frac{1}{n}]$ we have
  $h\ge 2^{s_n}\ge 2^k$, so that $\frac{h}{2^k}\ge 1$ on
  $[0,\frac{1}{n}]$, and it follows that $\frac{T^kh}{2^k}\ge 1$ on
  $[k\alpha,k\alpha+\frac{1}{n}]$.
  
  Now, given any positive integer $N$, we have $N\le
  s_{n-1}$ for some $n$. Then for each $k=s_{n-1}+1,\dots,s_n$ we have 
  $\frac{T^kh}{2^k}\ge 1$ on the interval 
  $[k\alpha,k\alpha+\frac{1}{n}]$. It follows that
  \begin{displaymath}
    \sum_{k=s_{n-1}+1}^{s_n}\tfrac{T^kh}{2^k}\ge 1
    \quad\mbox{ on the set }\quad
    \bigcup\limits_{k=s_{n-1}+1}^{s_n}[k\alpha,k\alpha+\tfrac{1}{n}]=
    s_{n-1}\alpha+\bigcup\limits_{k=1}^{M_n}[k\alpha,k\alpha+\tfrac{1}{n}]=
    [0,1],
  \end{displaymath}
  so that the series $\sum_{n=0}^\infty\frac{T^nh}{2^n}$ does not converge
  in measure.
\end{example}

\ssection{Other approaches to spectral theory}

There are different approaches to spectral theory of operators on
topological vector spaces, e.g.,~\cite{Waelbroeck:54}
and~\cite{Allan:65}. For example, Allan~\cite{Allan:65} defines the
spectrum of an element $x$ of locally-convex algebra $B$ as the set of
all $\lambda\in\mathbb C$ such that $\lambda e-x$ is not invertible or
the inverse is not bounded, where $y\in B$ is said to be bounded if
$\{c^n y^n\}_{n=1}^\infty$ is a bounded set for some real $c>0$. In
our terms, this means that $R_\lambda$ has finite spectral radius.
Allan's spectrum is, therefore, bigger than ours. Allan defines the
radius of boundedness of $\beta(x)$, which in our terms is exactly the
spectral radius, and he shows that $\beta(x)$ is less than or equal to
the geometrical radius of his spectrum. This result nicely complements
our Theorems~\ref{t:r-l}--\ref{t:r-nb} where we showed that a spectral
radius of an operator is greater that or equal to the geometrical
radius of the corresponding spectrum.

For example, if $X$ is a locally-convex space then it can be easily
verified that the collection of all continuous operators on $X$
equipped with the topology of uniform convergence on bounded sets is a
locally convex algebra.  For a base of convex neighborhoods of zero in
this algebra one can take the sets $\mathcal V_{\scriptscriptstyle
  A,U}$ of all continuous $T$ such that $T(A)\subseteq U$, where
$A\subseteq X$ is bounded and $U$ is a convex base zero neighborhood
in $X$. Therefore, the result of Allan is applicable in the setup of
Theorem~\ref{t:r-bb}.

\section{nb-bounded operators}   \lab{s:nb-bdd}

Since nb-boundedness is the strongest of the boundedness conditions we
have introduced, it is natural to expect that stronger results can be
obtained for nb-bounded operators.

\begin{numbered} \lab{n:nb-char}
  The following argument is often useful when dealing with nb-bounded
  operators. Suppose that $X$ and $Y$ are topological vector spaces
  and $T\colon X\to Y$ is nb-bounded, then $T(U)$ is a bounded set in
  $Y$ for some base zero neighborhood $U$. We claim that if $Y$ is
  Hausdorff, then $\bigcap_{n=1}^{\infty}\frac{1}{n}U\subseteq\Null T$.
  Indeed, it suffices to show that if $x\in\frac{1}{n}U$ for every
  $n\ge 1$ then $Tx$ belongs to every zero neighborhood $V$ of $Y$.
  But $T(U)\subseteq\alpha V$ for some positive $\alpha$
  (depending on $V$), and hence
  $Tx\in\frac{1}{n}T(U)\subseteq\frac{\alpha}{n}V\subseteq V$ whenever
  $n\ge\alpha$.
  
  It follows that if $T$ is one-to-one, then $U$ cannot contain any
  nontrivial linear subspaces. In particular, if $U$
  is convex then the locally bounded space $(X,U)$ is Hausdorff,
  hence quasinormable.
  In this case $T$ is a continuous operator from $(X,U)$ to $Y$, and,
  moreover, if $X=Y$, then $T$ is continuous as an operator from
  $(X,U)$ to $(X,U)$.
  
  In fact, many ``classical'' topological vector spaces have the
  property that every zero neighborhood contains a nontrivial linear
  subspace, e.g., topologies of pointwise or coor\-dinate-wise
  convergence, weak topologies, etc.
\end{numbered}

\begin{example} \lab{ex:ucc}
  {\it A topological vector space in which no base zero neighborhood
  contains a nontrivial linear subspace.} Let X be the space of all
  analytic functions on $\mathbb C$ equipped with the topology of
  uniform convergence on compact subsets of $\mathbb C$. The sets
  \begin{displaymath}
    U_{n,\varepsilon}=\bigl\{f\in X\mid \abs{f(z)}<\varepsilon
      \text{ whenever }\abs{z}\le n\bigr\} \quad 
      (n\ge 0\text{ and }\varepsilon>0) 
  \end{displaymath}
  form a zero neighborhood base of this topology. Clearly, no
  $U_{n,\varepsilon}$ contains a non-trivial linear subspace. Indeed,
  if there is a function $f$ in $X$ and a zero neighborhood
  $U_{n,\varepsilon}$ such that $\lambda f\in U_{n,\varepsilon}$ for
  every scalar $\lambda$, then $f(z)=0$ whenever $\abs{z}<n$, and it
  follows that $f$ is identically zero on $\mathbb C$. Note that this
  topology is generated by the countable sequence of seminorms
  $\norm{f}_n=\sup\limits_{\abs{z}\le n}\bigabs{f(z)}$; clearly
  $\norm{\cdot}_n$ is the Minkowski functional of $U_{n,1}$.
\end{example}

\begin{proposition} \lab{p:nb-not-invertible}
  If $X$ is a complete locally convex space then $X$ is locally
  bounded if and only if $X$ admits an nb-bounded bijection.
\end{proposition}

\begin{proof}
  If $X$ is locally bounded then the identity map is an nb-bounded
  bijection. Suppose that $T$ is an nb-bounded bijection on $X$. Then
  there exists a closed base zero neighborhood $U$ in $X$ such that
  $T(U)$ is bounded. Let $A=\overline{T(U)}$, then $A$ is convex,
  bounded, balanced, and absorbing. It follows that the space $(X,A)$
  is a locally convex and locally bounded, denote it by $X_A$. Notice
  also that the topology of $X_A$ is finer than the original topology
  on $X$ because $A$ is bounded. In particular, $X_A$ is Hausdorff.
  
  We claim that $X_A$ is complete. Indeed, if $(x_n)$ is a Cauchy
  sequence in $X_A$, then it is also Cauchy in the original topology
  of $X$, which is complete, so that $x_n$ converges to some $x$. Fix
  $\varepsilon>0$, then there exists $n_0$ such that
  $x_n-x_m\in\varepsilon A$ whenever $n,m\ge n_0$. Let $m\to\infty$,
  since $A$ is closed we have $x_n-x\in\varepsilon A$, i.e., $x_n\to
  x$ in $X_A$. Thus, $X_A$ is complete, hence Banach.
  
  Since $A$ is bounded, we can find $m$ such that $A\subseteq mU$.
  Then $T(A)\subseteq T(mU)\subseteq mA$, so that $T$ is bounded in
  $X_A$. Then $T^{-1}$ is also bounded in $X_A$ by the Banach Theorem,
  so that $U=T^{-1}\bigl(T(U)\bigr)\subseteq T^{-1}(A)\subseteq nA$
  for some $n>0$, hence $U$ is bounded.
\end{proof}

\begin{proposition}
  Let $T\colon X\to Y$ be an nb-bounded operator between Hausdorff
  topological vector spaces such that $X$ is not locally bounded. If
  \begin{enumerate}
    \item every zero neighborhood in $X$ contains a non-trivial linear
    subspace, or
    \item both $X$ and $Y$ are Fr\'echet spaces,
  \end{enumerate}
  then $T$ is not a bijection.
\end{proposition}

\begin{proof}
  If every zero neighborhood of $X$ contains a non-trivial linear
  subspace, then $T$ cannot be one-to-one by~\ref{n:nb-char}.  Suppose
  now that $X$ and $Y$ are Fr\'echet and assume that $T$ is a
  bijection. Let $S\colon Y\to X$ be the linear inverse of $T$. The
  Open Mapping Theorem implies that $S$ is continuous and hence
  bb-bounded. It follows that the identity operator of $X$ is
  nb-bounded being the composition of the nb-bounded operator $T$ and
  the bb-bounded operator $S$. But the identity operator is nb-bounded
  if and only if the space is locally bounded, a contradiction.
\end{proof}

\ssection{Weak topologies}

We are going to show that every operator which is nb-bounded relative
to a weak topology has to be of finite rank. In order to prove this we
need the following well-known lemma. For completeness we provide a
simple proof of it.

\begin{lemma} \lab{l:finiterank}
  Let $T$ be a linear operator on a vector space $L$, and let
  $f_1,\dots,f_n$ be linear functionals on $L$ such that $Tx=0$
  whenever $f_i(x)=0$ for every $i=1,\dots,n$. Then $T$ is a finite
  rank operator of rank at most $n$.
\end{lemma}

\begin{proof}
  Define a linear map $\pi$ from $L$ to ${\mathbb R}^n$ via
  $\pi(x)=(f_1(x),\dots,f_n(x))$.  Then the dimension of the range
  $\pi(L)$ is at most $n$.  Define also a linear map $\varphi$ from
  $\pi(L)$ to $L$ via $\varphi(\pi(x))=Tx$. It can be easily verified
  that $\varphi$ is well-defined. Then the range of $T$ coincides with
  the range $\varphi(\pi(L))$, which is of dimension at most $n$.
\end{proof}

\begin{proposition} \lab{p:nb-fntrank}
  Let $X$ be a locally convex space, and $T$ an operator on $X$ such
  that $T$ is nb-bounded with respect to the weak topology of $X$.
  Then $T$ is of finite rank.
\end{proposition}

\begin{proof}
  Suppose $T$ maps some weak base zero neighborhood $U=\bigl\{x\in X
  \mid\abs{f_i(x)}<1,\linebreak i=1,\dots,n\bigr\}$ ($f_1,\dots,f_n\in X'$),
  to a weakly bounded set. Since the weak topology is Hausdorff, it
  follows from~\ref{n:nb-char} that
  $\bigcap_{n=1}^{\infty}\frac{1}{n}U\subseteq\ker T$. In particular,
  $Tx=0$ whenever $f_i(x)=0$ for every $i=1,\dots,n$. Then
  Lemma~\ref{l:finiterank} implies that $T$ is a finite rank operator.
\end{proof}

\ssection{Spectra and spectral radii of nb-bounded operators}

\begin{proposition}  \lab{p:nb:spectra:equal}
  If $T$ is an nb-bounded operator on a topological vector space then
  $\sigma^{bb}(T)=\sigma^c(T)=\sigma^{nn}(T)=\sigma^{nb}(T)$.
\end{proposition}

\begin{proof}
  If $X$ is locally bounded then the result is trivial
  by~\ref{n:lb:spectra:equal}. Suppose that $X$ is not locally
  bounded, then, in view of~\ref{n:spectra-order}, it suffices to show
  that $\rho^{bb}(T)\subseteq\rho^{nb}(T)$. Let
  $\lambda\in\rho^{bb}(T)$, then $R_\lambda$ is bb-bounded. If
  $\lambda\neq 0$, then it follows from $R_\lambda(\lambda I-T)=I$
  that $R_\lambda=\frac{1}{\lambda}R_\lambda T+\frac{1}{\lambda}I$.
  Thus, $R_\lambda$ is a sum of an nb-bounded operator and a multiple
  of the identity operator, which yields $\lambda\in\rho^{nb}(T)$.  To
  finish the proof, it suffices to show that $\lambda=0$ necessarily
  belongs to $\sigma^{bb}(T)$ (and, therefore, to $\sigma^c(T)$,
  $\sigma^{nn}(T)$, and $\sigma^{nb}(T)$).  Indeed, if the resolvent
  $R_\lambda=T^{-1}$ were bb-bounded, then $I=T^{-1}T$ would be
  nb-bounded, which is impossible in a non-locally bounded space, a
  contradiction.
\end{proof}

\begin{proposition}  \lab{p:nb:radii:equal}
  If $T$ is an nb-bounded operator on a topological vector space, then
  $r_{bb}(T)=r_c(T)=r_{nn}(T)=r_{nb}(T)$.
\end{proposition}

\begin{proof}
  By Proposition~\ref{p:radii-order} it suffices to show that
  $r_{bb}(T)\ge r_{nb}(T)$. Since $T$ is nb-bounded, then $T(U)$ is a
  bounded set for some zero neighborhood $U$. Let $\nu>r_{bb}(T)$ and
  fix a zero neighborhood $V$. Then $\nu V$ is again a zero
  neighborhood. In particular, since the sequence $\tnun$ converges to
  zero uniformly on bounded sets, we have
  $\tnun\bigl(T(U)\bigr)\subseteq\nu V$ for all sufficiently large
  $n$. Then $\frac{T^{n+1}}{\nu^{n+1}}(U)\subseteq V$, so that $\tnun$
  converges to zero uniformly on $U$. Therefore $\nu\ge r_{nb}(T)$, so
  that $r_{bb}(T)\ge r_{nb}(T)$.
\end{proof}

\begin{numbered}  \lab{n:nb-notation}
  In view of Propositions~\ref{p:nb:spectra:equal}
  and~\ref{p:nb:radii:equal} we can write $\sigma(T)$\lab{a:sp-nb}
  instead of $\sigma^{bb}(T)$, $\sigma^c(T)$, $\sigma^{nn}(T)$, and
  $\sigma^{nb}(T)$ and $r(T)$\lab{a:sr-nb} instead of $r_{bb}(T)$,
  $r_c(T)$, $r_{nn}(T)$, and $r_{nb}(T)$.
\end{numbered}

We have established in Theorems~\ref{t:r-l}--\ref{t:r-nb} that under
some conditions the spectral radii of a linear operator are upper
bounds for the geometrical radii of the corresponding spectra. Of
course we would like to know when the equalities hold. It is well
known that the equality $\bigabs{\sigma(T)}=r(T)$ holds for every
continuous operator on a Banach space. Moreover, it was shown
in~\cite{Gramsch:66} that this equality also holds for every
continuous operator on a quasi-Banach space (a complete quasinormed
space).
Further, by means of Proposition~\ref{p:rc-convergence} the main result
of~\cite{Garibay:98} is equivalent to the following statement: {\it
  $r(T)=\bigabs{\sigma(T)}$ for every nb-bounded operator $T$ on a
  complete locally convex space\/}. Here we
present a direct proof of this. Our proof is a simplified version of
the proof of~\cite{Garibay:98}.
 
\begin{theorem}  \lab{t:rs=sr}
  If $T$ is an nb-bounded linear operator on a sequentially complete
  locally convex space, then $\bigabs{\sigma(T)}=r(T)$.
\end{theorem}

\begin{proof}
  Suppose $T(U)$ is bounded for some base zero neighborhood $U$. It
  follows from Propositions~\ref{p:nb:spectra:equal},
  \ref{p:nb:radii:equal}, and \ref{p:radii-order}, and
  Theorem~\ref{t:r-nb} that it suffices to show that
  $\bigabs{\sigma^{nn}(T)}\ge r_{nb}(T)$.  We are going to show that
  $T$ induces a continuous operator $\widetilde{T}$ on some Banach
  space such that
  $\sigma\bigl(\widetilde{T}\bigr)\subseteq\sigma^{nn}(T)\cup\{0\}$
  while $r\bigl(\widetilde{T}\bigr)\ge r_{nb}(T)$, and then appeal to
  the fact that the spectral radius of a continuous operator on a
  Banach space equals the geometrical radius of the spectrum.
  
  Consider $T$ as a continuous operator on the locally bounded space
  $X_\su=(X,U)$. Then $\sigma_\su(T)$ is defined
  by~\ref{n:lb:spectra:equal} and $r_\su(T)$ is defined
  by~\ref{n:lb:radii:equal}. We claim that $r_\su(T)\ge r_{nb}(T)$.
  To see this, suppose $r_\su(T)<\nu$, then
  $\frac{T^n}{\nu^n}(U)\subseteq U$ for all sufficiently large $n$.
  Let $V$ be a base zero neighborhood, then $T(U)\subseteq\alpha V$
  for some $\alpha>0$, so that
  $\frac{T^n}{\nu^n}(U)=\frac{T}{\nu}\frac{T^{n-1}}{\nu^{n-1}}(U)\subseteq
  \frac{1}{\nu}T(U)\subseteq\frac{\alpha}{\nu}V$ for sufficiently
  large $n$. This implies that $\nu\ge r_{nb}(T)$, and it follows that
  $r_\su(T)\ge r_{nb}(T)$.
  
  On the other hand, we claim that
  $\sigma_\su(T)\subseteq\sigma^{nn}(T)$. Suppose
  $\lambda\in\rho^{nn}(T)$, then $R_\lambda$ is nn-bounded with
  respect to some base ${\mathcal N}_0$ of zero neighborhood. We can
  assume without loss of generality that $U\in{\mathcal N}_0$, so that
  $R_\lambda(U)\subseteq\beta U$ for some $\beta>0$. It follows that
  $\lambda\in\rho_\su(T)$.
  
  Since $U$ is convex, the the space $X_\su$ is, in fact, a seminormed
  space. We can assume without loss of generality that it is a normed
  space, because otherwise we can consider 
  the quotient space $X_\su/(\Null T)$ and the quotient operator
  $\widehat{T}$ on this quotient space instead of $T$. Indeed, since   
  $\bigcap_{n=1}^{\infty}\frac{1}{n}U\subseteq\Null T$
  by~\ref{n:nb-char}, we conclude that the quotient space
  $X_\su/(\Null T)$ is Hausdorff. It follows then that $X_\su/(\Null
  T)$ is a normed space, and $\widehat{T}$ is norm bounded.
  The spectrum $\sigma_\su(T)$
  becomes even smaller when we substitute $T$ with $\widehat{T}$.
  Indeed, suppose $\lambda\in\rho_\su(T)$, then the resolvent
  $R_\lambda$ exists in $X_\su$ and is continuous. If $x\in\ker T$,
  then $x=R_\lambda(\lambda I-T)x=\lambda R_\lambda x$, so that
  $R_\lambda$ leaves $\ker T$ invariant, and, therefore, induces a
  quotient operator $\widehat{R_\lambda}$ on $X_\su/\ker T$ via
  $\widehat{R_\lambda}([x])=[R_\lambda x]$. Clearly,
  $\widehat{R_\lambda}$ is continuous: if $[x_n]\to [x]$ in
  $X_\su/\ker T$ then $x_n-z_n\to x$ in $X_\su$ for some
  $(z_n)_{n=1}^\infty$ in $\ker T$, so that $[R_\lambda
  x_n]=[R_\lambda (x_n-z_n)]\to[R_\lambda x]$. On the other hand,
  $r_\su(\widehat{T})\ge r_\su(T)$, because if
  $\nu>r_\su(\widehat{T})$ then
  $\frac{\widehat{T}^n}{\nu^n}([U])\subseteq[U]$ for all sufficiently
  large $n$, then $\tnun(U)\subseteq U+\ker T$, so that
  $\frac{T^{n+1}}{\nu^{n+1}}(U)\subseteq\frac{1}{\nu}T(U)\subseteq
  \frac{\alpha}{\nu}U$ for some $\alpha>0$. It follows that $\nu\ge
  r_\su(T)$ and, therefore, $r_\su(\widehat{T})\ge r_\su(T)$.
  
  Finally, we consider the completion $\widetilde{X}_\su$ of $X_\su$, and
  extend $T$ to a continuous linear operator $\widetilde{T}$ on the
  completion. The spectrum of $\widetilde{T}$ is smaller that the spectrum
  of $T$, because if $\lambda\in\rho_\su(T)$ then the resolvent
  $R_\lambda$ can be extended by continuity to $\widetilde{R_\lambda}$ on
  $\widetilde{X}$, and $\widetilde{R_\lambda}$ is a continuous inverse to
  $\lambda I-\widetilde{T}$, so that $\lambda\in\rho(\widetilde{T})$.  On the
  other hand, $r(\widetilde{T})\ge r_\su(T)$ because if
  $\nu>r(\widetilde{T})$ then
  $\frac{\widetilde{T}^n}{\nu^n}(\widetilde{U})\subseteq\widetilde{U}$ 
  for all
  sufficiently large $n$, which implies $\tnun(U)\subseteq U$ since
  $T$ is a restriction of $\widetilde{T}$ on $X$.
\end{proof}

\section{Compact operators}  \lab{s:compact}

As with bounded operators, there is more than one way to define
compact operators on an arbitrary topological vector space. A subset
of a topological vector space is called \term{precompact} if its
closure is compact. Given a linear operator $T$ on a topological
vector space, $T$ is called \term{Montel} if it maps every bounded set
into a precompact set and \term{compact\/} if it maps some
neighborhood into a precompact set.  To be consistent, we should have
probably called these operators ``b-compact'' and ``n-compact''
respectively, but the names ``Montel'' and ``compact'' are commonly
accepted.  Obviously, every compact operator is Montel and nb-bounded
(hence continuous); every Montel operator is bb-bounded.

\begin{numbered} \lab{r:r-compact}
  If $T$ is compact or Montel, then sequential completeness is
  not needed in Theorems~\ref{t:r-l}--\ref{t:r-nb}. Indeed, we used
  sequential completeness just once, namely, in the proof of
  Theorem~\ref{t:r-l} to justify the convergence of the sequence
  $R_{\lambda,n}x=\frac{1}{\lambda}\sum_{i=0}^n
  \frac{T^ix}{\lambda^i}$. But since the sequence $(R_{\lambda,n}x)_n$
  is Cauchy and, therefore, bounded, the sequence
  $(TR_{\lambda,n}x)_n$ has a convergent subsequence whenever $T$ is
  compact or Montel. Furthermore, it follows from
  $R_{\lambda,n+1}x=\frac{1}{\lambda}(I+TR_{\lambda,n})x$ that
  $(R_{\lambda,n}x)_n$ has a convergent subsequence hence converges.
\end{numbered}
  
Let $K$ be a compact operator on an arbitrary topological vector
space, and let $\sigma(K)$ and $r(K)$ be as in~\ref{n:nb-notation}.
It was proved in~\cite{Pech:91} that $\sigma(K)=\{0\}$ implies
$r_l(K)=0$.  In the following theorem we use the technique of
\cite{Pech:91} to improve this result by showing that in general
$r(K)\le\bigabs{\sigma(K)}$.

\begin{theorem}
  If $K$ is a compact operator on a Hausdorff
\marg{We don't need $K$ be injective, do we?}
  topological vector space $X$, then $r(K)\le\bigabs{\sigma(K)}$.
\end{theorem}

\begin{proof}
  Assume that
  $\bigabs{\sigma(K)}<r(K)$. Without loss of generality (by scaling
  $K$) we can assume that $\bigabs{\sigma(K)}<1<r(K)$.  Since $K$ is
  compact, there is a closed base
  zero neighborhood $U$ such that $\overline{K(U)}$ is compact. In
  particular $\overline{K(U)}$ is bounded, so that
  $\overline{K(U)}\subseteq \eta U$ for some $\eta>0$. We can assume
  without loss of generality that $\eta>1$.
  We define the following subsets of $U$:
  \[
    U_1=\overline{K(U)}\cap U,\quad
    U_{n+1}=K(U_n)\cap U \quad(n=1,2,\dots),\quad\mbox{and}\quad
    U_0=\bigcap\limits_{n=1}^\infty U_n.
  \]
  Notice, that $U_1$ is compact because $\overline{K(U)}$ is compact
  and $U$ is closed. Also, if $U_n$ is compact, then $K(U_n)$ is
  compact as the image of a compact set under a continuous operator.
  Therefore, every $U_n$ for $n\ge 1$ is compact. Using induction, we
  can show that the sequence $(U_n)$ is decreasing. Indeed,
  $U_1\subseteq U$ by definition, $U_2=K(U_1)\cap U\subseteq K(U)\cap
  U\subseteq U_1$, and if $U_n\subseteq U_{n-1}$, then
  $U_{n+1}=K(U_n)\cap U\subseteq K(U_{n-1})\cap U=U_n$. It follows
  also that $U_0$ is compact and contains zero.
  
  Notice that $K$ maps every balanced set to a balanced set. Since $U$
  is balanced, $U_n$ is balanced for each $n\ge 0$. If $A$ is a
  balanced subset of $U$, then obviously $A\subseteq(\eta A)\cap U$,
  and when we apply the same reasoning to $\frac{1}{\eta}K(A)$ instead
  of $A$ (which is also a balanced subset of $U$), we get
  $\frac{1}{\eta}K(A)\subseteq K(A)\cap U$. We use this to show by
  induction that $\frac{1}{\eta^n}K^n(U)\subseteq U_n$ for every $n\ge
  1$. Indeed, for $n=1$ we have $\frac{1}{\eta}K(U)\subset K(U)\cap
  U\subseteq U_1$.  Suppose $\frac{1}{\eta^n}K^n(U)\subseteq U_n$ for
  some $n\ge 1$, then
  \[
    \tfrac{1}{\eta^{n+1}}K^{n+1}(U)\subseteq\tfrac{1}{\eta}K(U_n)
    \subseteq K(U_n)\cap U=U_{n+1},
  \]
  which proves the induction step.
  
  Next, we claim that there exists an open zero neighborhood $V$ and an
  increasing sequence of positive integers $(n_j)$ such that
  $U_{n_j}\setminus V$ is nonempty for every $j\ge 1$. Assume for the
  sake of contradiction that for every open zero neighborhood $V$ we have
  $U_n\subseteq V$ for all sufficiently large $n$. Since
  $\frac{1}{2}U$ contains an open zero neighborhood, then there exists a positive
  integer $N$ such that $U_n\subseteq\frac{1}{2}U$ whenever $n\ge N$.
  This implies that $U_{N+m}=K^m(U_N)$ for all $m\ge 0$. Indeed, this
  holds trivially for $m=0$. Suppose that $U_{N+m}=K^m(U_N)$ for some $m\ge
  0$. Then $U_{N+m+1}=K(U_{N+m})\cap U=K^{m+1}(U_N)\cap U$, and this
  implies that $U_{N+m+1}=K^{m+1}(U_N)$ because
  $U_{N+m+1}\subseteq\frac{1}{2}U$. Now take any open zero neighborhood
  $V$, then $\frac{1}{\eta^N}V$ is again a zero neighborhood, and by
  assumption there exists a positive integer $M$ such that
  $U_n\subseteq\frac{1}{\eta^N}V$ whenever $n\ge M$. Let
  $n\ge\max\{M,N\}$, then
  \[
    V\supseteq\eta^N U_n=\eta^N K^{n-N}(U_N)\supseteq
    \eta^N K^{n-N}\left(\tfrac{1}{\eta^N}K^N(U)\right)=K^n(U),
  \]
  which contradicts the hypothesis $r_{nb}(K)=r(K)>1$.

  It follows from $U_{n_j}\setminus V\neq\emptyset$ for every $j\ge 1$
  that $U_n\setminus V\neq\emptyset$ for all sufficiently large $n$
  because $U_n$ is a decreasing sequence. Since $U_n\setminus V$
  is a decreasing sequence of nonempty compact sets, then
  $U_0\setminus V=\bigcap_{n=1}^\infty (U_n\setminus
  V)\neq\emptyset$, so that $U_0\neq\{0\}$.
  
  For every $n\ge 1$ we have $U_0\subseteq U_n$, it follows that
  $K(U_0)\subseteq K(U_n)$ and, therefore,
  $K(U_0)\subseteq\bigcap_{n=1}^\infty K(U_n)$.  Actually, the reverse
  inclusion also holds. To see this, let $y\in\bigcap_{n=1}^\infty
  K(U_n)$. Then $y=Kx_n$, where $x_n\in U_n\subseteq U_1$. Since $U_1$
  is compact, the sequence $(x_n)$ has a cluster point, i.e.,
  $x_{n_j}\to x$ for some subsequence $(x_{n_j})$ and some $x$. Since $K$ is
  continuous we have $y=Kx$. On the other hand, since every $U_{n_j}$ is
  closed we have $x\in U_{n_j}$, so that $x\in\bigcap_{n=1}^\infty
  U_{n_j}=U_0$. Thus $K(U_0)=\bigcap_{n=1}^\infty K(U_n)$.

  Next, we claim that $U_0\subseteq K(U_0)\subseteq\eta U_0$. Indeed,
  \[
    U_0=\bigcap\limits_{n=2}^\infty U_n=
    \bigcap\limits_{n=2}^\infty\bigl[K(U_{n-1})\cap U\bigr]
    \subseteq\bigcap\limits_{n=2}^\infty K(U_{n-1})=K(U_0).
  \]
  On the other hand, since $U_n$ are decreasing and $\eta>1$, we have
  $K(U_n)\subseteq K(U_{n-1})\subseteq\eta K(U_{n-1})$ and
  $K(U_n)\subseteq K(U)\subseteq\eta U$, so that $K(U_n)\subseteq\eta
  K(U_{n-1})\cap\eta U=\eta U_n$, and this implies $K(U_0)\subseteq
  K(U_n)\subseteq\eta U_n$ for every $n$. Thus $K(U_0)\subseteq\eta
  U_0$.
  
    Since $\overline{K(U)}$ is compact, hence bounded, then
  $\overline{K(U)}+\overline{K(U)}$ is also bounded.  Then there is a
  positive constant $\gamma$ such that
  $\overline{K(U)}+\overline{K(U)}\subseteq\gamma U$. Without loss of
  generality we can assume $\gamma\ge2$. It follows that
  \[
    U_1+U_1=\overline{K(U)}\cap U+\overline{K(U)}\cap U\subseteq
    \overline{K(U)}+\overline{K(U)}\subseteq\gamma U.
  \]
  We use induction to show that
  $U_n+U_n\subseteq\gamma U_{n-1}$. Indeed, since
  $A\cap B+C\cap D\subseteq(A+C)\cap(B+D)$ for any four sets
  $A$, $B$, $C$, and $D$, then
  \begin{multline*}
    U_{n+1}+U_{n+1}=K(U_n)\cap U_n+K(U_n)\cap U_n\\
    \subseteq\bigl[K(U_n)+K(U_n)\bigr]\cap(U_n+U_n)
    \subseteq K(U_n+U_n)\cap(U_n+U_n)\\
    \subseteq K(\gamma U_{n-1})\cap \gamma U_{n-1}
    =\gamma\bigl[K(U_{n-1})\cap U_{n-1}\bigr]=\gamma U_n.
  \end{multline*}
  Finally, $U_0+U_0\subseteq\bigcap_{n=1}^\infty (U_n+U_n)
  \subseteq\bigcap_{n=1}^\infty \gamma U_n=\gamma U_0$.
  
  Next, consider the set $F=\bigcup_{n=1}^\infty nU_0$. This set is
  closed under multiplication by a scalar, and $U_0+U_0\subseteq\gamma
  U_0$ implies that $F$ is a linear subspace of $X$. We consider the
  locally bounded topological vector space $(F,U_0)$ with multiples of
  $U_0$ as the base of zero neighborhoods. Since $U_0$ is balanced by
  definition, this topology is linear, and it is Hausdorff because
  $U_0$ is compact.  Also, it is finer than the topology on $F$
  inherited from $X$ because $U_0$ is compact and, therefore, bounded in
  $X$. 

  We claim that $(F,U_0)$ is complete. Indeed, if $(x_n)$ is a Cauchy
  sequence in $(F,U_0)$ then there exists $k>0$ such that $x_n\in
  kU_0$ for each $n>0$. Since $U_0$ is compact, the sequence $(x_n)$
  has a subsequence which converges to some $x\in kU_0$ in the
  topology of $X$. Moreover, $\lim\limits_{n\to\infty}x_n=x$ because
  the sequence $(x_n)$ is Cauchy in $X$. Fix $\varepsilon>0$, then
  there exists $n_0$ such that $x_n-x_m\in\varepsilon U_0$ whenever
  $n,m\ge n_0$. Let $m\to\infty$, since $U_0$ is is closed we have
  $x_n-x\in\varepsilon U_0$, i.e., $x_n\to x$ in $(F,U_0)$. Thus, 
  $(F,U_0)$ is complete and, therefore, quasi-Banach.
  
  It follows from $U_0\subseteq K(U_0)\subseteq\eta U_0$ that $F$ is
  invariant under $K$ and the restriction $\widetilde{K}=K|_F$ is
  continuous.  We claim that
  $\sigma(\widetilde{K})\subseteq\sigma(K)\cup\{0\}$.  Suppose that
  $\lambda\in\rho(K)$ and $\lambda\neq 0$, then $(\lambda I-K)$ is a
  homeomorphism, so that $(\lambda I-K)(U)$ is a closed zero
  neighborhood, and $\alpha U_1\subseteq(\lambda I-K)(U)$ for some
  positive real $\alpha$ because $U_1$ is bounded. Further, $\alpha
  K(U_1)\subseteq K(\lambda I-K)(U)\subseteq(\lambda I-K)K(U)$.
  Therefore
  \[
    \alpha U_2\subseteq \alpha K(U_1)\cap \alpha U_1
    \subseteq(\lambda I-K)K(U)\cap(\lambda I-K)(U),
  \]
  and since $\lambda I-K$ is one-to-one we get $\alpha
  U_2\subseteq(\lambda I-K)(K(U)\cap U)\subseteq (\lambda I-K)(U_1)$.
  Similarly, we obtain $\alpha U_{n+1}\subseteq(\lambda I-K)(U_n)$
  for all $n\ge 1$, and then $\alpha U_0\subseteq(\lambda
  I-K)(U_0)$. This implies that the restriction of $\lambda I-K$ to
  $F$ is onto, invertible, and the inverse is continuous. Thus,
  $\lambda\in\rho(\widetilde{K})$.
    
  In particular this implies that $\bigabs{\sigma(\widetilde{K})}\le
  \bigabs{\sigma(K)}<1$.  On the other hand, it follows from $U_0\subseteq
  K(U_0)$ that $U_0\subseteq\widetilde{K}^n(U_0)$ for all $n\ge 0$, so
  that $\widetilde{K}^n$ does not converge to zero uniformly on $U_0$,
  whence $r(\widetilde{K})=r_{bb}(\widetilde{K})\ge 1$. This produces a
  contradiction because it was proved in~\cite{Gramsch:66} that the
  spectral radius of a continuous operator on a quasi-Banach space
  equals the radius of the spectrum.
\end{proof}

\begin{corollary}
  If $K$ is a compact operator on a locally convex (or
  pseudo-convex) space, then
  $r(K)=\bigabs{\sigma(K)}$.
\end{corollary}

\section{Closed operators}

In certain situation one has to deal with unbounded linear operators
in Banach spaces. For example, the generator of a strongly continuous
operator semigroup is generally a closed operator with dense domain
(see e.g.~\cite{Hille:57,Dunford:58}). Through this section $T$ will
be a closed operator on a Banach space $X$ with domain $\mathcal
D(T)$. As usually, we define $\mathcal D(T^{n+1})=\{x\in\mathcal
D(T^n)\mid T^nx\in\mathcal D(T)\}$ and $D=\bigcap_{n=0}^\infty\mathcal
D(T^n)$. In case when $T$ is the infinitesimal generator of an
operator semigroup, $D$ is dense in the range of the semigroup,
which is usually assumed to be all of $X$. The set $D$ with the
locally-convex topology $\tau$ given by the sequence of norms
$\norm{x}_n=\sum_{k=0}^n\norm{T^kx}$ is a Fr\'echet space. Clearly,
$D$ is invariant under $T$, and the restriction operator $T_{|D}$ is
continuous because $x_\alpha\xrightarrow{\tau}0$ in $D$ implies
$\norm{Tx_\alpha}_n\le\norm{x_\alpha}_{n+1}\to 0$ for each $n$.

We investigate the relation between the spectral properties of the
original operator $T$ on $X$ and of the restriction $T_{|D}$ on $D$. A
different approach to this question can be found in~\cite{Wrobel:99}.
Recall that $\lambda\in\rho(T)$ if $R(\lambda;T)=(\lambda
I-T)^{-1}\colon X\to\mathcal D(T)$ exists (it is automatically bounded
by~\cite[Theorem~2.16.3]{Hille:57}), and $\sigma(T)=\mathbb
C\setminus\rho(T)$.

\begin{lemma} \lab{l:res}
  If $\lambda\in\rho(T)$ then $R(\lambda;T)$ is a bijection of $D$
  commuting with $T$ on $\mathcal D(T)$. Further, for each $n\ge 0$
  there is a constant $C_n$ such that $\norm{R(\lambda;T)x}_n\le
  C_n\norm{x}_{n-1}$ for each $x\in D$.
\end{lemma}

\begin{proof}
  It can be easily verified that $R(\lambda;T)$ is a bijection from
  $\mathcal D(T^n)$ onto $\mathcal D(T^{n+1})$ and, therefore, the
  restriction $R(\lambda;T)_{|D}$ is a bijection. Notice that since
  $R(\lambda;T)Tx=x$ for each $x\in\mathcal D(T)$ and
  $TR(\lambda;T)x=x$ for each $x\in X$, then
  \begin{eqnarray*}
    TR(\lambda;T)x & = & \lambda R(\lambda;T)x-x
          \mbox{ for each }x\in\mathcal D(T)\mbox{ and}\\
    R(\lambda;T)Tx & = & \lambda R(\lambda;T)x-x
          \mbox{ for each }x\in X,
  \end{eqnarray*}
  so that $T$ and $R(\lambda;T)$ commute on $\mathcal D(T)$. It also
  follows that for each $x\in D$ we have
  \begin{eqnarray*}
    T^2R(\lambda;T)x & = & \lambda TR(\lambda;T)x-Tx
                      = \lambda^2 R(\lambda;T)-\lambda x-Tx,\\
    T^3R(\lambda;T)x & = & \lambda^2 TR(\lambda;T)-\lambda Tx-T^2x
                           \lambda^3 R(\lambda;T)-\lambda^2 x-\lambda Tx-T^2x,\\
    &\vdots&\\
    T^kR(\lambda;T)x & = &\lambda^kR(\lambda;T)x-\lambda^{k-1}x
      -\lambda^{k-2}Tx-\dots-\lambda T^{k-2}x-T^{k-1}x.
  \end{eqnarray*}
  It follows that
  \begin{displaymath}
    \label{eq:a}
    \norm{T^kR(\lambda;T)x}\le\abs{\lambda}^k\norm{R(\lambda;T)x}+
    \abs{\lambda}^{k-1}\norm{x}+\abs{\lambda}^{k-2}\norm{Tx}+
    \dots+\norm{T^{k-1}x}
  \end{displaymath}
  for each $x\in D$, so that
  \begin{multline*}
    \norm{R(\lambda;T)x}_n=\sum_{k=0}^n\norm{T^kR(\lambda;T)x}\le\\
    \mu_n\norm{R(\lambda;T)x}+\mu_{n-1}\norm{x}+\mu_{n-2}\norm{Tx}+\dots
    +\mu_0\norm{T^{n-1}x},
  \end{multline*}
  where $\mu_k=1+\abs{\lambda}+\dots+\abs{\lambda}^k$. Since
  $\norm{R(\lambda;T)x}\le\norm{R(\lambda;T)}\norm{x}$ it follows that
  $\norm{R(\lambda;T)x}_n\le
  C_n\bigl(\norm{x}+\norm{Tx}+\dots+\norm{T^{n-1}x}\bigr)$ for some
  $C_n$, so that
  $\norm{R(\lambda;T)x}_n\le C_n\norm{x}_{n-1}$.
\end{proof}

\begin{proposition}
  The inclusion $\rho(T)\subseteq\rho^{nn}(T_{|D})$ holds.
  Moreover, if $D$ is dense in $X$ and $T$ is the smallest closed
  extension of $T_{|D}$, then $\rho(T)=\rho^{nn}(T_{|D})$.
\end{proposition}

\begin{proof}
  Suppose that $\lambda\in\rho(T)$ and consider the resolvent operator
  $R(\lambda;T)$ on $X$. Then $\norm{R(\lambda;T)x}_n\le
  C_n\norm{x}_{n-1}\le C_n\norm{x}_n$ hence $R(\lambda;T)_{|D}$ is
  nn-bounded and $\lambda\in\rho^{nn}(T_{|D})$.

  Suppose now that $D$ is dense in $X$, $T$ is the smallest closed
  extension of $T_{|D}$, and $\lambda\in\rho^{nn}(T_{|D})$. Then there
  exists an nn-bounded operator $R(\lambda;T_{|D})\colon D\to D$ such
  that $R(\lambda;T_{|D})(\lambda I-T)x=(\lambda
  I-T)R(\lambda;T_{|D})x=x$ for each $x\in D$. Then there is a
  constant $C>0$ such that
  $\norm{R(\lambda;T_{|D})x}=\norm{R(\lambda;T_{|D})x}_0\le
  C\norm{x}_0=C\norm{x}$ for each $x\in D$. It follows that
  $R(\lambda;T_{|D})$ can be extended to a bounded operator $R$ on
  $X$. Fix $x\in X$ and pick $(x_n)$ in $D$ such that $x_n\to x$. Then
  $R(\lambda;T_{|D})x_n\to Rx$ and $(\lambda
  I-T)R(\lambda;T_{|D})x_n=x_n\to x$. Since $\lambda I-T$ is closed we
  have $(\lambda I-T)Rx=x$. It follows, in particular, that $\lambda
  I-T$ is onto.
  
  Since $\norm{(\lambda I-T)x}\ge\frac{1}{C}x$ for each $x\in D$, it
  follows that for every nonzero $y\in X$ the pair $(y,0)$ doesn't
  belong to the closure of the graph of $\lambda I-T_{|D}$. But the
  closure of the graph of $\lambda I-T_{|D}$ is the graph of $\lambda
  I-T$ because $\lambda I-T$ is the smallest closed extension of
  $\lambda I-T_{|D}$. It follows that $\lambda I-T$ is one-to-one,
  hence $\lambda\in\rho(T)$.
\end{proof}

Suppose now that $S$ is a bounded operator on $X$ such that $D$ is
invariant under $S$ and $STx=TSx$ for each $x\in D$. Then
\begin{displaymath}
  \norm{Sx}_m=\sum_{k=0}^m\norm{T^kSx}
  \le\norm{S}\sum_{k=0}^m\norm{T^kx}=\norm{S}\cdot\norm{x}_m,
\end{displaymath}
so that $\norm{S_{|D}}_n\le\norm{S}$. Moreover, if $m\le k$ then
$\norm{x}_m\le\norm{x}_k$, so that the mixed seminorm
$\msn{k}{m}(S_{|D})\le\norm{S}$. It also follows from
Proposition~\ref{p:sr:lcs} that $r_{nn}(S_{|D})\le r(S)$.

Further, we claim that if $R=R(\lambda;T)$ for some
$\lambda\in\rho(T)$, then $r_{nb}(R_{|D})\le r(R)$. Indeed, recursive
application of Lemma~\ref{l:res} yields $\norm{R^nx}_k\le
M_k\norm{R^{n-k}x}$ for each $x\in D$ and $k\ge n$, where
$M_k=\Pi_{i=1}^kC_i$. It follows that the mixed seminorm
\begin{multline*}
  \msn{m}{k}(R_{|D}^n)=
    \sup\bigl\{\norm{R^nx}_k\mid x\in D,\,\norm{x}_m\le 1\bigr\}\\
    \le\sup\bigl\{M_k\norm{R^{n-k}x}\mid\norm{x}\le 1\bigr\}=
    M_k\norm{R^{n-k}}.
\end{multline*}
Therefore
$\lim_n\sqrt[n]{\msn{m}{k}(R_{|D}^n)}\le\lim_n\sqrt[n]{\norm{R^n}}
=r(T)$ for any $m,k\ge 0$.  Now Proposition~\ref{p:sr:lcs} yields
$r_{nb}(R_{|D})\le r(R)$.

%\bibliographystyle{alpha}
%\bibliography{tv}

\end{document}